\documentclass[11pt, reqno]{amsart}

\usepackage[T1]{fontenc}
\usepackage{lmodern}
\usepackage[a4paper,margin=28mm]{geometry}
\usepackage{amsmath,amssymb,amsthm,mathtools}
\usepackage{microtype}
\usepackage[colorlinks=true,linkcolor=blue,citecolor=blue,urlcolor=blue]{hyperref}

\numberwithin{equation}{section}
\newtheorem{theorem}{Theorem}[section]
\newtheorem{proposition}[theorem]{Proposition}
\newtheorem{lemma}[theorem]{Lemma}
\newtheorem{corollary}[theorem]{Corollary}
\newtheorem{example}[theorem]{Example}
\theoremstyle{definition}
\newtheorem{definition}[theorem]{Definition}
\newtheorem{remark}[theorem]{Remark}
\newtheorem{question}{Question}
\newtheorem{conjecture}{Conjecture}
\theoremstyle{remark}

\newcommand{\R}{\mathbb R}
\newcommand{\C}{\mathbb C}
\newcommand{\N}{\mathbb N}
\newcommand{\SE}{S_E}
\newcommand{\norm}[1]{\left\lVert #1\right\rVert}
\newcommand{\abs}[1]{\left\lvert #1\right\rvert}
\DeclareMathOperator{\Ree}{Re}

\newcommand{\floor}[1]{\left\lfloor#1\right\rfloor}
\newcommand{\dist}{\operatorname{dist}}

\newcommand{\Iso}{\operatorname{Iso}}

\begin{document}
\setcounter{page}{1}

\title[Step-isometries on  Banach spaces]{Step-isometries on  Banach spaces}

\author{Ruidong Wang}
\address{School of Mathematical Sciences, Tianjin University of Technology, Tianjin 300384, China}
\email{wangruidong@tjut.edu.cn}

\author{Li Tang}
\address{School of Mathematical Sciences, Tianjin University of Technology, Tianjin 300384, China}

\email{tangli@stud.tjut.edu.cn}

\author{Jianglong Luo}
\address{School of Mathematical Sciences, Tianjin University of Technology, Tianjin 300384, China}

\email{2100056943@qq.com}

\let\thefootnote\relax\footnote{ The study is supported by the Natural Science Foundation of China (Grant No.12271402) and Natural Science Foundation of Tianjin Municipal Science and Technology Commission (Grant No. 23JCYBJC00560, 22JCYBJC00420 )}

\date{Received: xxxxxx; Revised: yyyyyy; Accepted: zzzzzz.
}

\begin{abstract}
Step-isometries are maps preserving the integer part of pairwise distances on metric spaces. They arise naturally in the study of infinite random geometric graphs, where graph isomorphisms between LARG-type random graphs reduce to step-isometries on the underlying vertex sets. Several fundamental open problems concerning the rigidity, extendability, and continuity of step-isometries were raised in recent literature, notably by \v{S}emrl, and by Balogh, Walters, and Zs\'ak. In the present paper we resolve a number of these long-standing open questions.

We establish rigidity results for step-isometries defined on norm-dense subsets of both real and complex Banach spaces. For complex Banach spaces, we prove that every almost-surjective step-isometry on a norm-dense subset admits a unique extension to a surjective real-affine isometry of the whole space. In the real setting, we investigate three key properties $A$, $B$, $C$ formulated by \v{S}emrl, which ask whether certain classes of step-isometries extend to surjective affine isometries. We verify these properties for uniformly non-square Banach spaces, for spaces with a dense set of Fr\'echet-smooth points, and for strictly convex Banach spaces with separable dual. In particular, we confirm Conjecture 2 due to Balogh, Walters, and Zs\'ak without countability assumptions on the dense domain set. We further show that every bijective step-isometry defined on the entire Banach space is a homeomorphism, answering their Question 3 affirmatively. Specializing to $\ell_1^n$, we obtain that $\ell_1^n$ satisfies Property $A$ for all integers $n\ge 3$, settling an open problem posed by \v{S}emrl.

Turning to random geometric-graph applications, we settle Conjecture 1 of Balogh, Walters, and Zs\'ak: in a strictly convex Banach space with separable dual, typical countable dense subsets are strongly non-Rado. This conclusion holds in two standard senses of typicality: in the sense of Baire category for sequences in $X^{\mathbb N}$, and for independent samples drawn from full-support Borel measures that assign zero mass to every norm sphere. For such typical sets, independent random-graph realisations are almost surely non-isomorphic, and the random graphs are almost surely rigid with trivial automorphism group. Our work unifies and extends existing rigidity theory for step-isometries across finite- and infinite-dimensional Banach-space settings.

\end{abstract}
\maketitle

\section{Introduction}

Classical random-graph theory originates from the pioneering work of  Erd\H{o}s, R\'enyi, and Gilbert \cite{ER1960}. In  1959,    Erd\H{o}s and R\'enyi investigated the model \(G(n,m)\), the uniform random graph on $n$ vertices containing exactly $m$ edges \cite{ER1959}, while Gilbert proposed \(G(n,p)\), where each pair of distinct vertices is independently connected by an edge with probability \(p\in(0,1)\) \cite{Gilbert1959}. The latter is widely referred to in the literature as the Erd\H{o}s-R\'enyi random graph. In 1961, Gilbert constructed a geometric random-graph model whose vertices form a homogeneous planar Poisson point process; distinct vertices \(x,y\) are adjacent if and only if their Euclidean distance satisfies \(\|x-y\|_2<r\) for a fixed radius \(r>0\)  \cite{Gilbert,Penrose2003,Walters2011}.

The Erd\H{o}s-R\'enyi framework randomizes edge configurations while keeping the vertex set fixed. By contrast, Gilbert's geometric model randomizes vertex locations and defines adjacency via pairwise spatial distance. If one imposes the independent-edge mechanism on the countably infinite vertex set \(\mathbb{N}\), where each pair of natural numbers is connected independently   with probability
\(p\)), one obtains   the infinite random graph \(G(\mathbb{N},p)\).
In 1963, Erd\H{o}s--R\'enyi established that for every \(p \in (0,1)\), the random graph \(G(\mathbb{N},p)\) almost surely belongs to a single isomorphism class:
\[
\mathbb{P}\bigl(G(\mathbb{N},p)\cong R\bigr)=1 \qquad p\in(0,1),
\]
where the countably infinite graph \(R\) is independent of \(p\) \cite{ER1963,Cameron1997,Cameron2001}. Here \(\cong\) denotes graph isomorphism. Subsequently, Rado gave an explicit deterministic construction of \(R\) \cite{Rado}, which is now canonically known as the Rado graph. Recall that the Rado graph \cite{Rado, Cameron1997} is the unique countably infinite graph on vertex set \(\mathbb{N}\) satisfying the following property: for any two disjoint finite sets \(U, V \subset \mathbb{N}\)(where either may be empty), there exists a vertex \(w \in \mathbb{N} \setminus (U \cup V)\) adjacent to all vertices in \(U\) and to no vertex in \(V\) (any two graphs satisfying this property are isomorphic).

Motivated by this almost-sure isomorphism uniqueness, Bonato and Janssen combined geometric proximity with independent edge selection and asked when this uniqueness up to isomorphism persists under these geometric constraints \cite{Bonato2011, Janssen2010}.

We recall their model in the formulation adopted in the present paper \cite{Bonato2011}. Let $(M,d)$ be a metric space, where $d(x,y)$ denotes
the distance between $x$ and $y$.
Fix a countable vertex set $A\subseteq M$, a distance threshold $\delta>0$, and a probability parameter $p\in(0,1)$. The model is denoted by \(\operatorname{LARG}(A,\delta,p).\)
For a graph $G$ generated by this model, write $E(G)$ for its edge set.
Distinct vertices $x,y\in A$ satisfy
\[
\mathbb P\bigl(\{x,y\}\in E(G)\bigr)=p
\qquad\text{if }d(x,y)<\delta,
\]
and
\[
\mathbb P\bigl(\{x,y\}\in E(G)\bigr)=0
\qquad\text{if }d(x,y)\ge\delta,
\]
with independent realisations over all admissible pairs.
Here $\{x,y\}\in E(G)$ means that $x$ and $y$ are
adjacent, and $\mathbb P$ stands for probability.
Imposing geometric thresholds breaks the full symmetry enjoyed by \(G(\mathbb{N},p)\).
A central question therefore arises: are two independently sampled \(\operatorname{LARG}\)-graphs almost surely isomorphic despite these geometric constraints?

To analyse graph isomorphisms between such geometric random graphs, one must relate graph-theoretic path distance to the underlying metric distance. Graph distance, the minimal number of edges along a path connecting two vertices, is preserved under graph isomorphisms \cite{West2001}.

Graph distance is discrete, counting the number of edges (``steps''), in a connecting path, whereas the metric distance on the ambient space is continuous and measures actual spatial separation. In the \(\operatorname{LARG}\) model an edge appears precisely when metric distance lies strictly below \(\delta\). By the triangle inequality, if two vertices have graph distance \(k\), their metric distance is strictly less than \(k\delta\). since the metric length of any path with k edges cannot exceed the sum of edge-length upper bounds \(k\delta\). 
Consequently, any graph isomorphism, which preserves graph distances, must map pairs of points whose metric distance belongs to \([m\delta,(m+1)\delta)\) onto pairs lying in the same interval. Maps respecting the integer parts of continuous distances in this fashion are called  \emph{step-isometry}.

Bonato and Janssen introduced step-isometries to formalise such spatial constraints
\cite[Section~2]{Bonato2011}.
  Using this notion, Bonato and Janssen showed that the almost-sure isomorphism problem for these random graphs reduces to a problem concerning the existence of step-isometries \cite[Section~3]{Bonato2011}.
\begin{lemma}\label{yinli1}
Let \(A\) be a countable dense subset of a normed space. Then there exists a family \(\mathcal{H}\) comprising almost all graphs in \(\mathcal{G} = \mathcal{G}(A, p)\) such that whenever \(G_1, G_2 \in \mathcal{H}\) and a map \(f : A \to A\) induces an isomorphism between \(G_1\) and \(G_2\), the map \(f\) must be a step-isometry on \(A\).
\end{lemma}

As observed in \cite{Balogh}, Lemma \ref{yinli1}   does not   hold in in full generality as originally stated; this subtlety is not required for our purposes, and we refer the reader to \cite{Balogh}. for further discussion. An immediate consequence of Lemma \ref{yinli1} reads as follows.
\begin{corollary}\label{tuilun1}
Any countable dense subset of a normed space admitting at most countably many step-isometries is strongly non-Rado.
\end{corollary}

Within the Bonato-Janssen framework, isomorphism behaviour depends strongly on the geometry of the ambient metric space. In $\ell_\infty^n$, i.e., $\mathbb R^n$ equipped with the maximum norm
\(
\|x\|_\infty=\max_{1\le i\le n}|x_i|,
\) they established almost-sure isomorphism uniqueness
for countable dense sets satisfying suitable
non-integrality conditions.
By contrast, over the Euclidean plane,   independent
samples defined on any fixed countable dense set are almost surely non-isomorphic \cite{Bonato2011}. Their subsequent work \cite{Bonato2012,BonatoMetrics,Bonato2016} established almost
surely nonisomorphic behaviour for the hexagonal norm and excluded isomorphism behaviour in normed planes whose unit balls are strictly convex or polygonal, apart from parallelograms.

The general definition of step-isometry simplifies substantially when specialised to normed spaces with unit threshold, which is the setting of the present paper.
Let $X$ be a  normed space and let $A\subseteq X$
be a fixed countable dense subset.
Define $\widehat G(X,A)$ as the graph on vertex-set $A$ where
distinct vertices $x,y$ are adjacent precisely when
\(
\|x-y\|<1.
\)
Retaining each edge of $\widehat G(X,A)$ independently with probability
$p\in(0,1)$ yields the random graph $G_p(X,A)$,
and we write $\mathcal G_p(X,A)$ also for the corresponding
probability space.

\begin{definition}\label{budengjudingyi}
Let \(A\) be a dense subset of normed space \(X\). A bijection \(f:A\to A\), \(f\) is called a \textit{step-isometry} if
 $$\lfloor \|x - y\| \rfloor = \lfloor \|f(x) - f(y)\| \rfloor$$
 for all \(x, y \in A\).

Equivalently, \(f\) is a step-isometry provided that for any \(k \in \mathbb{N}\), one hase \(\|x - y\| < k\) if and only if \(\|f(x) - f(y)\| < k\).
\end{definition}
From the Definition \ref{budengjudingyi} one readily obtains the estimate
\begin{equation}\label{eq:one-error}
\bigl|\|f(x)-f(y)\|-\|x-y\|\bigr|<1.
\end{equation}
A step-isometry therefore preserves the integer part of distances without necessarily preserving distances exactly. Determining which normed-space geometries force every step-isometry to be a genuine isometry becomes a central geometric problem for understanding isomorphism phenomena of random geometric graphs.

Balister et al. introduced the terminology
of \emph{Rado sets} and \emph{strongly non-Rado sets}  
to encode isomorphism behaviour purely as a property of the vertex set $A$ \cite{BBGLW}.
A countable dense set $A\subseteq X$ is  
\emph{Rado} if  for every $p\in(0,1)$, two independent
graphs $G_1,G_2$ sampled from $\mathcal G_p(X,A)$
satisfy
\(
\mathbb P(G_1\cong G_2)=1.
\)
The set $A$ is called \emph{strongly non-Rado} if  for every
$p\in(0,1)$,
\(
\mathbb P(G_1\cong G_2)=0.
\)
These concepts describe isomorphism statistics for random graphs built over a fixed vertex set; they do  not  require the resulting graph to be isomorphic to the abstract Rado graph R.

We emphasise an important conceptual distinction.
For a fixed vertex set $A$, statements involving ``almost surely'' refer to random edge selections within \(G_p(X,A)\). By contrast, ``almost all'' or ``typical'' vertex sets $A$ concern random or generic choice of $A$ itself, and require either a probability measure or a Baire-category formulation. Whenever needed we will state explicit additional conditions on $A$, such as absence of repeated pairwise distances.

Balister et al. gave a classification for typical countable dense subsets in finite-dimensional normed spaces \cite[Theorems~1]{BBGLW}.
\begin{theorem}\label{BBdingli}
\cite{BBGLW} Let \(X\) be a finite-dimensional normed space not linearly isometric to \(\ell_\infty^n\). Then for every \(0 < p < 1\), almost every countable dense subset \(A\) is strongly non-Rado.
\end{theorem}
Combined with the results of Bonato and Janssen, Theorem \ref{BBdingli} identifies $\ell_\infty^n$ as, up to linear isometry, the unique $n$-dimensional normed space for which almost all countable dense sets are Rado. Balister et al. also constructed examples of countable dense sets that
are neither Rado nor strongly non-Rado. Thus the behaviour of typical sets and the existence of exceptional sets represent separate classification problems.

Existing theory on step-isometries develops along two main directions, which we treat in separate subsections below. The first direction studies bijective step-isometries \(A\to A\), which are directly tied to graph isomorphisms and hence to the Rado / strongly non-Rado dichotomy. The second direction investigates step-isometries \(A\to X\) mapping into the ambient space; here the main themes are affine rigidity, extendability, continuity, and representation theorems. The present paper contributes to both lines of research.

\subsection{Bijective step-isometries \(A\to A\) and graph isomorphisms}

As mentioned above, Lemma \ref{yinli1} yields Corollary \ref{tuilun1}, which gives a sufficient condition for a countable dense set to be strongly non-Rado: it suffices that the ambient space admits only countably many bijective step-isometries on $A$. In several normed spaces there exist abundant step-isometries that fail to be isometries.

In several normed spaces there exist abundant step-isometries that fail to be isometries. For example, take any increasing continuous bijection \(g:[0,1]\to[0,1]\) satisfying \(g(0)=0\), \(g(1)=1\). One readily verifies that \[f(x)=\lfloor x\rfloor+g(x-\lfloor x\rfloor)\] defines a step-isometry on \(\mathbb{R}\). Maps of this form are referred to as standard step-isometries on \(\mathbb{R}\). On the other hand, in Euclidean spaces \(\mathbb{R}^n\) (\(n \ge 2\)), every step-isometry defined on a countable dense subset is an isometry \cite{Dilworth, SV}. This motivates the general problem: characterise normed spaces for which every bijective step-isometry on every countable dense subset is an isometry.

Any such isometry defined on a dense subset \(S\) admits extension to a bijective isometry on the whole ambient space \(X\), which is affine by the Mazur-Ulam Theorem (see e.g., \cite{MazurUlam}). If \(X\)  is finite-dimensional and possesses this rigidity property, then any dense subset \(S\) has at most countably many isometries. Appealing to Corollary \ref{tuilun1}, one deduces that almost all (or all) countable dense subsets of \(X\) are strongly non-Rado.

In this finite-dimensional setting Balister et al. \cite{BBGLW} confirmed that \(\ell_\infty^n\) is the sole space where almost all countable dense sets are Rado. In other words, whenever a finite-dimensional normed space is not isometric to \(\ell_\infty^n\),  almost every countable dense subset is strongly non-Rado.

Results in infinite dimensions are more varied. Bonato, Janssen and Quas \cite{Bonato2019, Bonato2021} proved that almost all countable dense sets are Rado in \(c\), \(c_0\) and \(C([0,1])\) (the spaces of convergent sequences, null sequences, and continuous functions on \([0,1]\), respectively, all equipped with the sup norm). Angel and Spinka \cite{Angel} studied the circle quotient \(\mathbb{R}/L\): for irrational \(L\) every countable dense set is strongly non-Rado, whereas for rational \(L\)  almost all countable dense sets are Rado.

More recently, Balogh, Walters and Zs\'ak \cite{Balogh} established the following rigidity theorem for reflexive Banach spaces
\begin{theorem}\label{dingli1}
Let \(X\) be a strictly convex reflexive Banach space with \(\dim X\ge2\), let \(A\subset X \) be  dense, and let \(f\) be a step-isometry on \(A\). Then \(f\) is an isometry.
\end{theorem}
\begin{corollary}\label{tuilun2}
Let \(X\) be a strictly convex reflexive Banach space with \(\dim X\ge2\), and let \(A\) be a countable dense set with no repeated pairwise distances. Then \(A\) is strongly non-Rado.
\end{corollary}

An interesting feature of Theorem \ref{dingli1} is that it holds without separability assumptions on $X$, unlike the graph-model context which inherently requires countability of the vertex-set $A$. Balogh et al. also constructed a Banach space in which every countable dense set is strongly non-Rado. Their argument proceeds by first extending a step-isometry from the dense domain to the whole space and then proving that global step-isometries are isometries.

In finite dimensions Balister et al. \cite{BBGLW} proved that each step-isometry on a dense set \(A\) extends uniquely to a step-isometry on the whole   space, relying on the norm compactness of the unit ball. Norm-compactness of the unit ball fails in infinite dimensions, \cite{Balogh} exhibit an explicit counterexample in \(c_0\) demonstrating that extensions need not exist. This motivates imposing additional structural hypotheses such as reflexivity, where weak compactness of bounded sets facilitates extension arguments. Nevertheless, since step-isometries are defined via the norm topology, infinite-dimensional extension arguments are substantially more involved than their finite-dimensional counterparts. Balogh, Walters and Zs\'ak \cite{Balogh} further proved the following extension theorem within reflexive Banach spaces.
\begin{theorem}
Let \(f\) be a step-isometry defined on a dense subset \(A\) of a reflexive Banach space \(X\). Then \(f\) admits a unique extension to a step-isometry \(\bar{f}\) on the whole space \(X\).
\end{theorem}
Existence of an extension does not automatically guarantee exact distance preservation: non-isometric global step-isometries already appear in finite-dimensional examples \cite{BBGLW}. Strict convexity changes the picture dramatically and permits carrying over finite-dimensional reasoning to the reflexive infinite-dimensional setting, as encapsulated in Theorem \ref{dingli1}.

 While Theorem \ref{dingli1} gives strong rigidity under strict convexity, the class of reflexive Banach spaces includes many non-strictly-convex examples. Balogh, Walters and Zs\'ak therefore formulated the following open question \cite{Balogh}.
\begin{question}\label{question1}
For which reflexive spaces \(X\) must every step-isometry be an  isometry?
\end{question}

The behaviour of typical countable dense vertex sets in the non-reflexive regime remains another major open direction. Balogh, Walters and Zs\'ak \cite{Balogh} showed that every infinite-dimensional Banach space with a \(1\)-symmetric basis contains at least one countable dense Rado set; this applies in particular to \(\ell_1\). Their construction yields existence of some Rado set, but leaves open the statistical behaviour of typical dense subsets  of \(\ell_1\).
\begin{question}\label{q:l1}
For a typical countable dense subset of \(\ell_1\), is the set Rado or strongly non-Rado?
\end{question}

To push the rigidity theory beyond reflexivity, Balogh, Walters, and Zs\'ak proposed two conjectures concerning strictly convex Banach spaces with separable dual \cite{Balogh}.
\begin{conjecture}\label{conj:typical}
Let \(X\) be a strictly convex Banach space with separable dual. Then `typical' countable dense subsets of \(X\) are \textit{strongly non-Rado}.
\end{conjecture}
\begin{conjecture}\label{conj2}
Let \(X\) be a strictly convex Banach space with separable dual, and  let \(f\) be a step-isometry on a countable dense set \(A\subset X\). Then \(f\)   extends to a step-isometry on all of \(X\).
\end{conjecture}

Counterexamples in \cite{Balogh} show that such extensions may fail in general, owing to discontinuity of the constructed map. Step-isometries are known to be continuous on many classical Banach spaces, yet the \(c_0\) example draws a clear distinction between step-isometries defined only on dense incomplete domains versus those defined globally on the complete Banach space. This observation prompted another problem posed in \cite{Balogh}.
\begin{question}\label{question3}
Suppose   \(X\) is a Banach space and   \(f\) is a step-isometry defined on the whole space \(X\). Must \(f\) be continuous?
\end{question}

Broadly speaking, one would like to exploit geometric properties of Banach spaces \cite{Clarkson,James,Huff,DGK,Asplund,BP} to characterize the existence of Rado or strongly non-Rado sets in infinite-dimensional spaces. This program also feeds back into random-graph theory. Step-isometry as a subject has generated a substantial body of results. 

\subsection{Step-isometries \(A\to X\) and affine rigidity}
\v{S}emrl \cite{Semrl} initiated the systematic study of step-isometries \(f : A \to X\)   on   dense subsets \(A\) of a Banach space \(X\) without assuming bijectivity of \(f\). He obtained representation theorems for dense-domain step-isometries on real \(\ell_p\)-spaces, \(1<p\le\infty\).
\begin{theorem}
Let \(1<p<\infty\), and let \(X=\ell_p^n\), \(n\ge2\), or
\(X=\ell_p\). Suppose that \(A\subset X\) is dense and
\(f:A\to X\) is a step-isometry. If \(X\) is infinite-dimensional,
assume additionally that \(f\) is almost surjective (a map \(f : A \to X\) is said to be almost surjective if there exists a positive real number \(\delta\) such that for every \(y \in X\) there exists \(x \in A\) such that \(\|f(x) - y\| \le \delta\)). Then
\[
    f(x)=Ux+b \qquad x\in A,
\]
where \(b\in X\) and \(U:X\to X\) is a surjective linear isometry.
Moreover, \(U\) is orthogonal when \(p=2\), and a signed coordinate
permutation otherwise.
\end{theorem}
\begin{theorem}
Let \(X=\ell_\infty^n\), \(n\ge2\), or \(X=\ell_\infty\).
Suppose that \(A\subset X\) is dense and \(f:A\to X\)
is a step-isometry with dense image. Then
\[
    f(x)=\bigl(\varepsilon_k\tau_k(x_{\sigma(k)})\bigr)_k+b
    \qquad x\in A,
\]
where \(b\in X\), \(\sigma\) is a permutation of the coordinate
index set, \(\varepsilon_k\in\{-1,1\}\), and each
\(\tau_k:\mathbb R\to\mathbb R\) is a standard step-isometry.
\end{theorem}
For  finite-dimensional \(\ell_p\)-spaces. Theorem 1.8  improves results from Balister et al. \cite{BBGLW} by removing the surjectivity assumptions. These representation results raise natural follow-up questions: what affine rigidity persists if one drops the almost-surjectivity or dense-image hypotheses? Recall that in infinite dimensions a linear isometric embedding need not be surjective. From
\v{S}emrl's theorems \cite{Semrl} one deduces that every bijective step-isometry \(T:S\to S\) on a countable dense subset of \(\ell_p\), \(1<p<\infty\) must be an isometry, so such sets   are strongly non-Rado. Nevertheless, \v{S}emrl left the corresponding problems for \(L_p\)-spaces and the general classification of Banach spaces  rigidity properties.
\v{S}emrl \cite{Semrl} concluded his paper with the following open problems.
\begin{question}\label{question4}
 Characterize those Banach spaces \(X\) for which every step-isometry \(f:A\to X\)  with \(A\subset X\) dense  extends to a surjective affine isometry of \(X\).
\end{question}
\begin{question}\label{question5}
 The same conclusion should hold  for every almost surjective step-isometry \(f:A\to X\).
\end{question}
\begin{question}\label{question6}
 The same conclusion should hold  for every step-isometry \(f:A\to X\) with dense image.
\end{question}
We say that a Banach space \(X\) has \textbf{Property A}, if every step-isometry \(f:A\to X\)  with \(A\subset X\) dense  extends to a surjective affine isometry of \(X\). A Banach space \(X\) has \textbf{Property B}, if every almost surjective step-isometry \(f:A\to X\) extends to a surjective affine isometry of \(X\). A Banach space \(X\) has \textbf{Property C}, if every  step-isometry \(f:A\to X\) with dense image extends to a surjective affine isometry of \(X\).

Properties \(A\), \(B\) and \(C\) exactly to Problems 5.1--5.3 in \v{S}emrl \cite{Semrl}, and satisfy the implications
\begin{equation}\label{implications}
\mathrm Property\, A\Longrightarrow \mathrm Property\, B\Longrightarrow \mathrm Property\, C.
\end{equation}
An affine isometric embedding \(f(x)=Ux+b\) where \(U\) is merely an isometric (non-surjective) embedding does \textbf{not} satisfy the defining requirement for \(Property\, A\).

 \v{S}emrl \cite{Semrl} observes that \(\ell_1^2\), which is isometric to \(\ell_\infty^2\), fails \(\mathrm Property\, A\), and asked whether \(\mathrm Property\, A\) holds for \(\ell_1^n\) whenever \(n\ge3\).

In the present paper we study geometric conditions that guarantee affine rigidity, extendability, and continuity for step-isometries, and derive consequences for strong non-Rado behaviour. In particular, we address several of the open problems listed above. 
The organisation of the paper is as follows.

Section 2 collects necessary background material and auxiliary lemmas concerning step-isometries, subdifferentials of norms, and related geometric tools, which will be invoked throughout our subsequent arguments.

In Section 3, we study almost-surjective step-isometries defined on norm-dense subsets of complex Banach spaces. There we establish our rigidity theorem for the complex setting, showing that every such map extends uniquely to a surjective real-affine isometry.

Section 4 treats step-isometries on $\ell_p$ spaces for $1<p<\infty$. In particular, we settle an open problem of \v{S}emrl by proving that $\ell_1^n$ satisfies Property $A$ whenever $n\ge 3$.

Section 5 investigates the affine-rigidity properties $A$, $B$, and $C$ introduced by \v{S}emrl. We verify these properties for several classes of real Banach spaces, including uniformly non-square spaces and strictly convex spaces with separable dual. This section also confirms Conjecture~2 of Balogh, Walters, and Zs\'ak, without imposing countability on the underlying dense domain.

Section 6 addresses Question~3 raised by Balogh, Walters, and Zs\'ak. We prove that every bijective step-isometry defined on the whole Banach space is a homeomorphism, giving an affirmative answer to this open question.

Section 7 is devoted to random-geometric-graph applications. There we resolve Conjecture~1 of Balogh, Walters, and Zs\'ak: in strictly convex Banach spaces with separable dual, typical countable dense subsets are strongly non-Rado, both in the sense of Baire category and for suitable Borel sampling measures. Consequences for almost-sure rigidity and non-isomorphism of random geometric graphs are also derived.

Finally, we close with a brief concluding remark and outlook.

\section{Main Lemmas and Proposition of  Step-Isometries}
Let $E$ be a real or complex Banach space, and let $A$ be a norm-dense subset. Suppose that  \(f:A\to E\) is a  step-isometry. Define
$$\Omega_f=\{\varphi\in E^*: \varphi(f(x))\le \varphi(f(y)) \text{ for all } x,y\in A \text{ with } \phi(x)<\varphi(y)\}.$$


Recall that a map \(g:E\to F\) between  Banach spaces is an
\(\varepsilon\)-nearisometry if
\[
 \abs{\norm{g(x)-g(y)}-\norm{x-y}}\leq\varepsilon
 \qquad(x,y\in E).
\]

For a normed space $X$,  we write
\[
  B_{X}=\{x\in X:\|x\|\leq1\},
  \qquad
  S_{X}=\{x\in X:\|x\|=1\},
\]
for its closed unit ball and unit sphere, respectively.

We first work over the real scalar field. A real continuous linear
functional \(\varphi\) is a \emph{unit supporting functional} at
\(w\ne0\) if \(\norm{\varphi}=1\) and
\(\varphi(w)=\norm{w}\). Such a functional exists by Hahn--Banach theorem.
A unit vector \(u\) is said to be \emph{G\^ateaux smooth} if the norm is
G\^ateaux derivative at \(u\), in which case the derivative is given by a unique real continuous linear functional \(j(u)\). Equivalently, \(j(u)\) is the unique unit
supporting functional at \(u\). We write $\mathrm{s.m.}(E)$ for the collection of all unit vectors in $E$ at which the norm is G\^ateaux-smooth. For any $x\in X$, we write
$$\partial_x=\{f\in S_{X^*}: f(x)=\|x\|\}.$$

The following theorem, which will be used later, can be found in \cite{Holmes}.
\begin{theorem}[\cite{Holmes}]\label{lem:generic}
Let \(E\) be a nonzero separable real Banach space. The G\^ateaux smooth
points are dense in its unit sphere.
\end{theorem}

The following result has been proved in \cite{Dilworth} and \cite{SV}.

\begin{theorem}\label{thm:stability}
Let $h:X\to X$ be an $\varepsilon$-isometry with $h(0)=0$. There is a surjective linear isometry $W$ such that
\[\norm{h(x)-Wx}\le2\varepsilon\] for every $x$ if either $X$ is finite-dimensional, or $h$ is almost surjective.
\end{theorem}

\begin{lemma}\label{lem:normalization}
Let \(E\) be a real or complex Banach space, and let \(A\) be a norm-dense subset of $E$. Suppose that
 \(f:A\to E\) is a  step-isometry.
Fix \(a_0\in A\), set $A_0=A-a_0$, and define
\[
 f_0(x)=f(x+a_0)-f(a_0),\quad \forall x\in A_0.
\]
If either $X$ is finite-dimensional or $f$ is almost surjective. Then there exists a surjective real-linear isometry \(U:E\to E\) such that the
mapping $$F=U^{-1}f_0:A_0\to E$$ is a step-isometry  satisfying $F(0)=0$ and
\begin{equation}\label{eq:normalized}
  \norm{F(x)-x}\leq6,\quad\forall x\in A_0.
\end{equation}
\end{lemma}

\begin{proof}
It is clear that \(f_0\)
is a  step-isometry  and maps \(0\) to \(0\).

For each \(z\in E\), let $a_z=z$, if \(z\in A_0\). If \(z\not\in A_0\), choose  \(a_z\in A_0\) such that
\(\norm{a_z-z}<1\).
Define $$g(z)=f_0(a_z).$$ Then
\[
 g|_{A_0}=f_0,\qquad g(0)=0,\qquad g(E)=f_0(A_0).
\]
In particular, \(g\) is almost surjective when  \(f\) is almost surjective.

For any $z,w\in E$,  by
\eqref{eq:one-error},
\begin{align*}
 \abs{\norm{g(z)-g(w)}-\norm{z-w}}
 &\leq
 \abs{\norm{f_0(a_z)-f_0(a_w)}-\norm{a_z-a_w}}\\
 &\quad+\abs{\norm{a_z-a_w}-\norm{z-w}}\\
 &<1+\norm{a_z-z}+\norm{a_w-w}<3 .
\end{align*}

Apply Theorem~\ref{thm:stability} to the underlying real Banach spaces.
This yields a surjective real-linear isometry \(U\) satisfying
$$\norm{g(z)-Uz}\leq6, \quad \forall z\in E.$$ Therefore, for all  \(x\in A_0\),
$$\norm{F(x)-x}=\|U^{-1}f_0(x)-x\|=\|f_0(x)-U(x)\|=\|g(x)-U(x)\|\leq6$$ which establishes
\eqref{eq:normalized}. Since $f_0$ is a step-isometry and  \(U^{-1}\) is an isometry, it follows that $F$ is also a step-isometry.

\end{proof}

\begin{lemma}\label{lem:moving}
Let \(E\) be a real Banach space, \(u\in s.m.(E)\cap S_E\), and let  \(\{w_k\}_{k\in \mathbb N}\subset E\).
If  $\lim\limits_{k\to \infty}w_k= u,$ then
\[
 \lim\limits_{k\to \infty}\varphi_k(z)= \varphi(z).
\]
for every $z\in E$, $\varphi\in \partial_z$
and  every \(\varphi_k \in \partial_{w_k}\).

\end{lemma}

\begin{proof}
Fix \(z\in E\). For any $\epsilon>0$,   the subdifferential supporting inequalities yield
\[
 \frac{\norm{w_k}-\norm{w_k-tz}}{t}
 \leq\varphi_k(z)\leq
 \frac{\norm{w_k+tz}-\norm{w_k}}{t}, \quad \forall t>0.
\]

Since \(u\in s.m.(E)\cap S_E\), we may pick $\delta>0$ such that
\[
 \frac{\norm{u}-\norm{u-\delta z}}{\delta}
 \leq\varphi(z)\leq
 \frac{\norm{u}-\norm{u-\delta z}}{\delta} +\epsilon,
\]
and
\[
 \frac{\norm{u+\delta z}-\norm{u}}{\delta}-\epsilon
 \leq\varphi(z)\leq
 \frac{\norm{u+\delta z}-\norm{u}}{\delta}.
\]
As $w_k\to u$, there exists $N_0\in \mathbb N$ satisfying
 $$\big| \frac{\norm{w_k+\delta z}-\norm{w_k}}{\delta}- \frac{\norm{u+\delta z}-\norm{u}}{\delta}\big|<\epsilon,\quad \big| \frac{\norm{w_k}-\norm{w_k+\delta z}}{\delta}- \frac{\norm{u}-\norm{u+\delta z}}{\delta}\big|<\epsilon$$
for all $k\ge N_0$. Combining these estimates,
$$\varphi_k(z)-\varphi(z)\le \frac{\norm{w_k+tz}-\norm{w_k}}{t}-\frac{\norm{u+\delta z}-\norm{u}}{\delta}+\epsilon<2 \epsilon,$$
and
$$\varphi_k(z)-\varphi(z)\ge \frac{\norm{w_k}-\norm{w_k+\delta z}}{\delta}- \frac{\norm{u}-\norm{u+\delta z}}{\delta}-\epsilon>-2 \epsilon.$$

Since $\epsilon>0$ is arbitrary, we conclude \(
 \lim\limits_{k\to \infty}\varphi_k(z)= \varphi(z).
\)

\end{proof}

\begin{lemma}\label{lem:fixed}
Let \(E\) be a real Banach space, let \(u\in s.m.(E)\) and let $\{v_k\}_{k\in \mathbb N}\subset E$.
Suppose that \(t_k>0\), \(t_k\to\infty\), and \(v_k/t_k\to u\).
Then
\begin{equation}\label{eq:fixed}
 \lim\limits_{k\to \infty}(\norm{v_k-a}-\norm{v_k-b})=
  j(u)(b-a).
\end{equation}
for all  \(a,b\in E\).
\end{lemma}

\begin{proof}
For arbitrary  \(a,b\in E\). Choose unit supporting functionals \(\varphi_k\in \partial_{v_k-a}\) and \(\phi_k\in \partial_{v_k-b}\), respectively. Observe that
$$\varphi_k\in \partial_{\frac{v_k-a}{t_k}}\quad and \quad \phi_k\in \partial_{\frac{v_k-b}{t_k}}.$$

Since
$$\lim\limits_{k\to \infty}\frac{v_k-a}{t_k}=\lim\limits_{k\to \infty}\frac{v_k-b}{t_k}=u,$$ and

\[
 \phi_k(b-a)\leq
 \norm{v_k-a}-\norm{v_k-b}
 \leq\varphi_k(b-a).
\]
By Lemma~\ref{lem:moving}, we conclude $\lim\limits_{k\to \infty}(\norm{v_k-a}-\norm{v_k-b})=
  j(u)(b-a)$.
\end{proof}

\begin{proposition}\label{prop:order}
Let \(E\) be a real Banach space, and let \(A\) be a norm-dense subset of $E$. Suppose that
\(F:A\to E\) is a step-isometry satisfying
\[
 \sup_{x\in A}\norm{F(x)-x}<\infty.
\]

Then, for every G\^ateaux smooth unit vector \(u\in S_E\) and all $x,y\in A$,  if $j(u)(x)<j(u)(y)$, it follows that
$$j(u)(F(x))\leq j(u)(F(y)).$$
\end{proposition}

\begin{proof}
Fix a G\^ateaux smooth unit vector \(u\in S_E\) and take \(x,y\in A\) satisfying \(j(u)(x)<j(u)(y)\). Choose
\(c\in\R\) such that \(j(u)(x)<c<j(u)(y)\).
For each   integer \(k\in \mathbb N\) with $k>|c|$, select \(z_k\in A\)
satisfying
\[
 \norm{z_k-(k+c)u}<1/k .
\]
Since $u\in s.m.(E)$, we have
\begin{align*}
\lim\limits_{k\to \infty}\big(\norm{z_k-x}-(k+c)+j(u)(x)\big)&=\lim\limits_{k\to \infty}\big(\norm{z_k-x}-(k+c)-\frac{\|u-\frac{1}{k+c}x\|-1}{\frac{1}{k+c}}\big)\\
&=\lim\limits_{k\to \infty}\big(\|z_k-x\|-\|(k+c) u-x\|\big)=0.
\end{align*}
The same reasoning applied to \(y\) yields
\[
 \lim\limits_{k\to \infty}\big(\norm{z_k-y}-(k+c)+j(u)(y)\big)=0.
\]
Accordingly, there exists $N_0\in \mathbb N$ such that
\[
 \norm{z_k-x}>k,\qquad \norm{z_k-y}<k,
\]
for all $k\ge N_0$.
Since $F$ is a step-isometry,
\[
 \norm{F(z_k)-F(x)}\geq k>
 \norm{F(z_k)-F(y)}.
\]
From the uniform bound $ \sup_{x\in A}\norm{F(x)-x}<\infty$, we deduce
$$\lim\limits_{k\to\infty}\big(\frac{F(z_k)}{k}- u\big)=\lim\limits_{k\to\infty}\big(\frac{F(z_k)-z_k}{k}+ \frac{z_k-(k+c)u}{k}+\frac{cu}{k}\big)=0.$$
Apply Lemma~\ref{lem:fixed} with \(v_k=F(z_k)\),
\(a=F(x)\), and \(b=F(y)\), we get
$$ j(u)(F(y))-j(u)(F(x))=\lim\limits_{k\to \infty}(\norm{F(z_k)-F(x)}-\norm{F(z_k)-F(y)})\ge 0
.$$
This completes the proof.

\end{proof}

\begin{lemma}\label{lem:continuous}
Let \(E\) be a real Banach space, and let \(A\) be a norm-dense subset of $E$. Suppose that
\(F:A\to E\) is a step-isometry satisfying
\[
 \sup_{x\in A}\norm{F(x)-x}=M<\infty.
\]
Let $\phi\in S_{X^*}$ satisfy
$\phi(F(x))\leq \phi(F(y))$ for all $x,y \in A$ with $\phi(x)<\phi(y)$. Define
\[
g_\phi(t)=\sup\{\phi(Fx):x\in A,\ \phi(x)<t\},
\]
and
$$D_\phi=\{t\in \R: g_\phi \text{ is continious at } t\}.$$
Then $g_\phi(t)$ is finite and nondecreasing, and
\begin{equation}\label{sandwich}
t-M\le g_\varphi(t)\le t+M,\qquad
g_\phi(\phi(x))\le\phi(Fx)\le g_\phi(\phi(x)+).
\end{equation}
Moreover, $\R\setminus D_\phi$ is countable. For every $r\in D$ and every sequence $\{x_n\}_{n\in \N}\subset A$ with $\lim\limits_{n\to\infty}\phi(x_n)=r$, we have
\begin{equation}\label{env-limit}
\lim\limits_{n\to\infty}\phi(Fx_n)= g_\phi(r).
\end{equation}
\end{lemma}
\begin{proof}
Set $$S_t = \{x \in A : \phi(x) < t\}.$$ Since \(\phi(A)\) is dense in \(\mathbb{R}\),  \(S_t \) is nonempty.

For any \(x \in S_t\), since \(\|Fx - x\| \le M\), \(\phi(x) < t\) and
\[
\phi(Fx) = \phi(x) + \phi(Fx - x) \le \phi(x) + |\phi(Fx - x)| \le \phi(x) + \|\phi\|\|Fx - x\|,
\]
it follows that \(g_\phi(t) \le t + M\).

To establish the lower bound, because \(\phi(A)\) is dense in \(\mathbb{R}\), we can choose a sequence \(\{x_n\} \subset A\) such that \(\phi(x_n) < t\) and \(\lim\limits_{n\to\infty}\phi(x_n) = t\).

For each \(n\),  \(g_\phi(t) \ge \phi(Fx_n)\), and
\[
\phi(Fx_n) = \phi(x_n) + \phi(Fx_n - x_n) \ge \phi(x_n) - |\phi(Fx_n - x_n)| \ge \phi(x_n) - M.
\]
Therefore, $$g_\phi(t) \ge \sup\phi(Fx_n)\ge t - M.$$

Clearly, \(g_\phi\) is   nondecreasing.

Now fix \(x \in A\). We prove
 $$g_\phi(\phi(x))\le\phi(Fx)\le g_\phi(\phi(x)+).$$
 If \(y \in S_{\phi(x)}\), then \(\phi(y) < \phi(x)\). By hypothesis, \(\phi(Fy) \le \phi(Fx)\). Taking the supremum over all such $y$, we get \(g_\phi(\phi(x)) \le \phi(Fx)\).

Conversely, pick any  \(s > \phi(x)\). By density, there exists \(z \in A\) satisfying $$\phi(x) < \phi(z) < s.$$ Then $\phi(Fx) \le \phi(Fz).$ Since \(\phi(z) < s\), we have \(z \in S_s\), so \(\phi(Fz) \le g_\phi(s)\). Hence $$\phi(Fx) \le g_\phi(s),\quad \forall s > \phi(x).$$ Taking the infimum as \(s \downarrow \phi(x)\), we obtain \(\phi(Fx) \le g_\phi(\phi(x)+)\).

Finally, let \(\{x_n\} \subset A\) satisfy \(\phi(x_n) \to r \in D_\phi\). For every \(n\),
\[
g_\phi(\phi(x_n)) \le \phi(Fx_n) \le g_\phi(\phi(x_n)+).
\]
Since \(g_\phi\) is continuous at $r$,
$$\lim\limits_{n\to\infty}g_\phi(\phi(x_n))=\lim\limits_{n\to\infty}g_\phi(\phi(x_n)+)=g_\phi(r),$$ so  \(\lim\limits_{n\to\infty}\phi(Fx_n)= g_\phi(r)\).

This completes the proof.

\end{proof}

\begin{lemma}\label{lem:finitecontinuous}
Let \(E\) be a real Banach space, and let \(A\) be a norm-dense subset of $E$. Suppose that
\(F:A\to E\) is a step-isometry satisfying
\[
 \sup_{x\in A}\norm{F(x)-x}=M<\infty.
\]
 Let $\psi_1,\ldots,\psi_m\in S_{E^*}$ be linearly independent with $m\ge2$, and let $a=\{a_i\}_{i=1}^m\in \R^m$ satisfy $a_i\not=0$ for all $1\le i\le m$. Set $\varphi=\sum_i a_i\psi_i$. If $$\{\varphi, \psi_1,\ldots,\psi_m\}\subset \Omega_F,$$ then
\[
\psi_i(Fx)=\psi_i(x)\quad(1\le i\le m),\qquad\text{ and } \qquad \varphi(Fx)=\varphi(x),
\]
for all $x\in A$.
\end{lemma}
\begin{proof}
Since $-\psi\in \Omega_F$ whenever $\psi\in \Omega_F$,
 we may assume that  $a_i>0$ for all \(1\le i\le m.\)

Define the map $g_{\psi_i}$ as in Lemma 2.5,  and let $D_i$ denote the set of continuity points of $g_{\psi_i}$. Set
\[
D=D_1\times\cdots\times D_m,\qquad
a=(a_1,\ldots,a_m),
\]
and define
\[
G(t)=\sum_{i=1}^m a_i g_{\psi_i}(t_i)
\qquad \forall t=(t_1,\ldots,t_m)\in\mathbb R^m.
\]

First, we claim that $G(s)\le G(t)$ for all $s,t\in D$ satisfying $a\cdot s<a\cdot t$.
To verify this claim, consider the linear map
\[
\Psi:X\longrightarrow\mathbb R^m,\qquad
\Psi(x)=(\psi_1(x),\ldots,\psi_m(x)).
\]
Since $\{\psi_1,\ldots,\psi_m\}$ is linearly independent, the map $\Psi$ is surjective. Therefore,
for any $s,t\in D$ with $a\cdot s<a\cdot t$, there exist $u,v\in X$ such that
\[
\Psi(u)=s,\qquad \Psi(v)=t.
\]
By density of   $A$, pick sequence $\{x_n\},\{y_n\}\subset A$ such that $\lim\limits_{n\to \infty}x_n= u$ and $\lim\limits_{n\to \infty}y_n= v$.
Then
\[
\lim\limits_{n\to \infty}\varphi(x_n)= a\cdot s,
\qquad
\lim\limits_{n\to \infty}\varphi(y_n)= a\cdot t.
\]
As $a\cdot s<a\cdot t$, there exists $N_0\in \N$ such that  $\varphi(x_n)<\varphi(y_n)$ for all  $n\ge N_0$. Since $\varphi\in \Omega_F$,
\[
\sum_{i=1}^m a_i\psi_i(Fx_n)
\le
\sum_{i=1}^m a_i\psi_i(Fy_n),
\]
holds for all  $n\ge N_0$.Applying Lemma \ref{lem:continuous}, we obtain
$$G(s)=\lim\limits_{n\to \infty}\sum_{i=1}^m a_i\psi_i(Fx_n)\le \lim\limits_{n\to \infty}\sum_{i=1}^m a_i\psi_i(Fy_n)=G(t).$$
Second, we  prove that
\begin{equation}\label{mixing-increments}
g_i(s)-g_i(r)=s-r,
\qquad \forall r,s\in D_i, \quad \forall 1\le i\le m.
\end{equation}
Fix an index $i$, choose $j\ne i$, and take $r,s\in D_i$ with $r<s$. Set
\[
d=s-r>0,\qquad \Delta=g_i(s)-g_i(r).
\]
For any $0<\delta<d$, define
\[
L_+=\frac{a_i}{a_j}(d+\delta)>0.
\]
The set \[
B=\bigcup_{n=0}^{\infty}(\mathbb R\setminus D_j-nL_+)
\]
is countable, so we may select $v_+\in \mathbb R\setminus B$.
It follows that
\[
v_++nL_+\in D_j\qquad n=0,1,2,\ldots.
\]

Fix arbitrary $c_\ell\in D_\ell$ for $\ell\ne i,j$. For every $n\ge0$,
choose $\alpha^{(n)},\beta^{(n)}\in D$ satisfying
\[
\begin{aligned}
\alpha_i^{(n)}&=s,
&
\alpha_j^{(n)}&=v_++nL_+,\\
\beta_i^{(n)}&=r,
&
\beta_j^{(n)}&=v_++(n+1)L_+,
\end{aligned}
\]
and $\alpha_\ell^{(n)}=\beta_\ell^{(n)}=c_\ell$. Then
\[
a\cdot\beta^{(n)}-a\cdot\alpha^{(n)}
=-a_i d+a_jL_+
=a_i\delta>0.
\]
Hence $G(\alpha^{(n)})\le G(\beta^{(n)})$, which yields
\[
a_i\Delta
\le
a_j\bigl[
g_{\psi_j}(v_++(n+1)L_+)-g_{\psi_j}(v_++nL_+)
\bigr]. \qquad n=0,1,2,\ldots.
\]
For any $N\in \N$, using(\ref{sandwich}) from Lemma \ref{lem:continuous},
\[
\begin{aligned}
\Delta
&\le \frac{a_j}{Na_i}
\sum_{n=0}^{N-1}
\bigl[g_j(v_++(n+1)L_+)-g_j(v_++nL_+)\bigr]
\\
&=\frac{a_j}{Na_i}
\bigl[g_j(v_++NL_+)-g_j(v_+)\bigr]\\
&\le\frac{a_j}{Na_i}(NL_++2M)
\\
&=d+\delta+\frac{2a_jM}{Na_i}.
\end{aligned}
\]
On the other hand, set
\[
L_-=\frac{a_i}{a_j}(d-\delta)>0,
\]
and choose $v_-\in\mathbb R$ such that
\[
v_-+nL_-\in D_j\qquad(n\ge0).
\]
With the remaining coordinates fixed as before, take
\[
\begin{aligned}
\alpha_i^{(n)}&=r,
&
\alpha_j^{(n)}&=v_-+(n+1)L_-,\\
\beta_i^{(n)}&=s,
&
\beta_j^{(n)}&=v_-+nL_-.
\end{aligned}
\]
Now
\[
a\cdot\beta^{(n)}-a\cdot\alpha^{(n)}
=a_i d-a_jL_-
=a_i\delta>0.
\]
It follows that
\[
a_j\bigl[
g_j(v_-+(n+1)L_-)-g_j(v_-+nL_-)
\bigr]
\le a_i\Delta.
\]

We now derive the lower bound,
\[
\begin{aligned}
\Delta
&\ge \frac{a_j}{Na_i}
\sum_{n=0}^{N-1}
\bigl[g_j(v_-+(n+1)L_-)-g_j(v_-+nL_-)\bigr]
\\
&=\frac{a_j}{Na_i}
\bigl[g_j(v_-+NL_-)-g_j(v_-)\bigr]
\\
&\ge\frac{a_j}{Na_i}(NL_--2M)
\\
&=d-\delta-\frac{2a_jM}{Na_i}
.
\end{aligned}
\]

Consequently,
\[
d-\delta-\frac{2a_jM}{Na_i}\le\Delta\le d+\delta+\frac{2a_jM}{Na_i}.
\]
Since $N\in \N$ and $\delta>0$ are arbitrary, we conclude that
\[
\Delta=d.
\]
This establishes \eqref{mixing-increments}.

It follows that $g_i(t)-t$ is constant on $D_i$. Since $D_i$ is dense in $\R$ and $g_i(t)$ is nondecreasing, there exists some $\beta_i\in \R$ such that
\[
g_i(t)=t+\beta_i\qquad(t\in\mathbb R).
\]
 In particular, $g_i$ is continuous on $\R$.

By Lemma \ref{lem:continuous},
\[
\psi_i(Fx)=g_i(\psi_i(x))
=\psi_i(x)+\beta_i
\qquad \forall x\in A.
\]
Setting $x=0$ and using $F(0)=0$, we obtain $\beta_i=0$. Therefore
\[
\psi_i(Fx)=\psi_i(x)
\qquad(1\le i\le m,\ x\in A).
\]
Finally,
\[
\varphi(Fx)
=\sum_{i=1}^m a_i\psi_i(Fx)
=\sum_{i=1}^m a_i\psi_i(x)
=\varphi(x).
\]
This completes the proof.
\end{proof}

\section{Almost Surjective Step-Isometries on Complex Banach Spaces}
In this section, we prove that every almost surjective step-isometry from a norm-dense
subset of a complex Banach space into that space extends uniquely to a
surjective real affine isometry.

\begin{lemma}\label{lem:plane}
Let \(E\) be a complex Banach space, let \(\Phi\in E^*\), let \(A\) be a norm-dense subset $E$ with
\(0\in A\). Suppose that   \(F:A\to E\) satisfies \(F(0)=0\).
If for every \(\theta\in\R\) and all \(x,y\in A\),  $\Ree(e^{-i\theta}\Phi(x))<\Ree(e^{-i\theta}\Phi(y))$ implies $\Ree(e^{-i\theta}\Phi(F(x)))\le\Ree(e^{-i\theta}\Phi(F(y)))$,
then there exists a real number \(a_\Phi\geq0\) such that
\[
 \Phi(F(x))=a_\Phi\Phi(x),\qquad\forall x\in A.
\]
\end{lemma}

\begin{proof}
For \(x,y\in A\), set
\[
 \alpha=\Phi(y-x),\qquad \beta=\Phi(F(y)-F(x)).
\]
When \(\alpha\ne0\), the hypothesis implies that  $\Ree(e^{-i\theta}\alpha)>0$ implies $\Ree(e^{-i\theta}w)\geq0$.
It follows that \(\beta=\lambda \alpha\) for some \(\lambda\geq0\).

We next prove that \(\Phi(x)=\Phi(y)\) implies
\(\Phi(F(x))=\Phi(F(y))\). Since \(\Phi(A)\) is dense in \(\C\),
we may choose \(z_1,z_2\in A\) such that
\(\Phi(z_1-x)\) and \(\Phi(z_2-x)\) are real-linearly independent.
Applying the preceding conclusion, there exist $\lambda_1, \lambda_2\ge 0$ satisfying
$$\Phi(F(z_1)-F(x))=\lambda_1 \Phi(z_1-x)\quad and \quad \Phi(F(z_1)-F(y))=\lambda_2 \Phi(z_1-y)=\lambda_2 \Phi(z_1-x).$$
Hence $$\Phi(F(x)-F(y))=(\lambda_2-\lambda_1) \Phi(z_1-x)\in \R\Phi(z_1-x).$$
By symmetry, we similarly obtain $$\Phi(F(x)-F(y))\in \R\Phi(z_2-x).$$

Therefore,
\[
 \Phi(F(x)-F(y))\in
 \R\Phi(z_1-x)\cap\R\Phi(z_2-x)=\{0\}.
\]
Thus  for any $x\in A$, put $t=\Phi(x)$. Define the map
\(h:\Phi(A)\to\C\)
  by
$$h(t)=\Phi\circ F(x).$$ This map is well-defined and satisfies \(h(0)=0\). Moreover, for distinct \(s, t\in \Phi(A)\), there exists $\lambda_{s,t}\ge 0$ such that
$$h(s)-h(t)=\lambda_{s,t}(s-t).$$

Fix two real-linearly independent elements \(s_1,s_2\in \Phi(A)\). Then there exist $\alpha_1, \alpha_2, \alpha\ge 0$ such that
$$h(s_1)=\alpha_1s_1, \quad h(s_2)=\alpha_2s_2$$
and
\[
 h(s_2)-h(s_1)=\alpha (s_2-s_1),\qquad\lambda\geq0.
\]
Real-linear independence forces
$$\alpha=\alpha_1=\alpha_2.$$
Set $a_\Phi=\alpha$. For any nonzero \(s\in \Phi(A)\), at least one of \(s_1,s_2\) is
real-linearly independent of \(s\).
Repeating the  argument for that pair gives \(h(s)=a_\Phi s\). Consequently,
\[
 \Phi(F(x))=a_\Phi\Phi(x),\qquad\forall x\in A.
\]
This completes the proof.

\end{proof}

For a complex Banach space \(E\), a \emph{complex supporting functional}
at \(u\in\SE\) is a complex-linear \(\Phi\in E^*\) with
\(\norm{\Phi}=1\) and \(\Phi(u)=1\).
Its real part is a real unit supporting functional.
Conversely, a real continuous functional \(\varphi\) is the real part of the
unique complex-linear functional
\[
 \Phi(x)=\varphi(x)-i\,\varphi(ix),
\]
and the two functional norms agree.
It follows that real G\^ateaux smoothness at \(u\) is equivalent to
uniqueness of its complex supporting functional.
We denote the latter by \(\Phi_u\).

\begin{theorem}\label{thm:criterion}
Let \(E\) be a complex Banach space, and let \(A\) be a norm-dense subset of $E$ with \(0\in A\). Suppose that
\(F:A\to E\) be a step-isometry satisfying
\[
  F(0)=0,\qquad \sup_{x\in A}\norm{F(x)-x}<\infty.
\] Assume further that there exists a family
\(\mathcal U\subset s.m.(E)\cap S_E\)  whose complex
supporting functionals separate points, i.e.
\[
 \bigcap_{u\in\mathcal U}\ker\Phi_u=\{0\}.\]
Then \(F(x)=x\) for all \(x\in A\).
\end{theorem}

\begin{proof}
Set $$M=\sup_{x\in A}\norm{F(x)-x}<\infty.$$

Fix \(u\in\mathcal U\). For every \(\theta\in\R\), the vector
\(e^{i\theta}u\) is G\^ateaux smooth, and its real supporting
functional is
\[
 j(e^{i\theta}u)=\Ree(e^{-i\theta}\Phi_u).
\]

Applying Proposition~\ref{prop:order}  to
each such rotated vector shows that the hypotheses of
Lemma~\ref{lem:plane} are satisfied. Consequently,
\[
 \Phi_u(F(x))=a_u\Phi_u(x),\qquad\forall x\in A,
\]
for some \(a_u\geq0\).
Since \(\norm{\Phi_u}=1\), we have
\[
 \abs{(a_u-1)\Phi_u(x)}
 =\abs{\Phi_u(F(x)-x)}\leq M, \qquad\forall x\in A.
\]
The set \(\Phi_u(A)\) is dense in \(\C\), in particular, unbounded. This forces \(a_u=1\).
Since this holds for every \(u\in\mathcal U\), the point-separating condition yields \(F(x)=x\) for all \(x\in A\).

\end{proof}

\begin{lemma}\label{lem:separating}
If \(E\) is a nonzero separable complex Banach space, the complex
supporting functionals at its smooth unit vectors separate points.
\end{lemma}

\begin{proof}
Given \(z\ne0\), choose, by Theorem~\ref{lem:generic}, a smooth unit vector
\(u\) such that \(\norm{u-z/\norm z}<1/2\). Since
\(\norm{\Phi_u}=1\) and \(\Phi_u(u)=1\),
\[
 \bigl|\Phi_u(z/\norm z)\bigr|
 \ge1-\norm{u-z/\norm z}>1/2.
\]
Thus \(\Phi_u(z)\ne0\).
\end{proof}

\begin{corollary}\label{cor:separable}
Let \(E\) be a separable complex Banach space. Let \(A\subset E\) be
norm-dense with \(0\in A\), and let \(F:A\to E\) be a step-isometry with
\[
 F(0)=0,\qquad \sup_{x\in A}\norm{F(x)-x}<\infty.
\]
Then \(F\) is the identity on \(A\).
\end{corollary}

\begin{proof}
For \(E\ne\{0\}\), combine Lemma~\ref{lem:separating} with
Theorem~\ref{thm:criterion}. The zero space is immediate.

\end{proof}

\begin{lemma}\label{lem:reduction}
Let \(E\) be a complex Banach space, let \(A\) be a norm-dense subset of $E$
with \(0\in A\), and let \(F:A\to E\) be a mapping satisfying \(F(0)=0\).
For each \(x_0\in A\), there exist a closed separable complex subspace
\(Y\subset E\) and a countable set \(D\subset A\cap Y\) such that
\[
 0,x_0\in D,\qquad \overline {D}^{\|\cdot\|}=Y,\qquad F(D)\subset Y.
\]

\end{lemma}

\begin{proof}
Set
\[
 Y_0=\overline{\operatorname{span}_{\C}}^{\|\cdot\|}\{x_0,F(x_0)\}.
\]
Inductively, suppose we are given a closed separable complex subspace \(Y_n\). Choose
a countable dense set \(Q_n\subset Y_n\). Since \(A\) is  norm-dense in $E$, for each \(q\in Q_n\) and each
positive integer \(m\), there exists \(a_{q,m}\in A\) such that
\(\norm{a_{q,m}-q}<1/m\). Define
\[
 D_n=\{a_{q,m}:q\in Q_n,\ m\ge1\},
\]
and
\[
 Y_{n+1}
 =\overline{\operatorname{span}_{\C}}^{\|\cdot\|}
    \bigl(Y_n\cup D_n\cup F(D_n)\bigr).
\]
Then \(Y_{n+1}\) is separable and \(Y_n\subset\overline{D_n}^{\|\cdot\|}\).
Now set
\[
 Y=\overline{\bigcup_{n\ge0}Y_n}^{\|\cdot\|},
 \qquad D=\{0,x_0\}\cup\bigcup_{n\ge0}D_n.
\]
Then \(Y\) is a closed separable
complex subspace of $E$. By construction, $$D\subset A\cap Y,\quad
F(D)\subset Y,$$ and \(Y_n\subset \overline {D}^{\|\cdot\|}\) for every $n\in \mathbb N$.
Consequently, \(\overline {D}^{\|\cdot\|}=Y\).

This completes the proof.

\end{proof}

\begin{proposition}\label{prop:allbounded}
Let \(E\) be any complex Banach space, let \(A\) be a norm-dense
subset of $E$ with \(0\in A\), and let \(F:A\to E\) be a step-isometry satisfying
\[
 F(0)=0,\qquad \sup_{x\in A}\norm{F(x)-x}<\infty.
\]
Then \(F(x)=x\) for all \(x\in A\).
\end{proposition}

\begin{proof}
Fix \(x_0\in A\), and apply  Lemma~\ref{lem:reduction} to obtain \(D,Y\).
The restriction $$F|_D:D\to Y$$ is a step-isometry on a norm-dense
subset of the separable complex Banach space \(Y\), maps 0 to 0,
and lies at bounded distance from the identity map.
Corollary~\ref{cor:separable} yields \(F(x_0)=x_0\). Since \(x_0\) is
arbitrary, the desired conclusion follows.

\end{proof}

\begin{theorem}
Let \(E\) be any complex Banach space, let \(A\) be a norm-dense subset of $E$,
and let \(f:A\to E\) be an almost surjective step-isometry.
Then there exist a surjective real-linear isometry \(U:E\to E\) and \(b\in E\)
such that
\[
 f(x)=Ux+b, \qquad\forall x\in A.
\]
In particular, \(f\) has a unique extension to a surjective real affine
isometry of \(E\).
\end{theorem}

\begin{proof}
Fix \(a_0\in A\). by
Lemma~\ref{lem:normalization}, there exists a surjective real-linear
isometry \(U:E\to E\) and a normalized map
\[
 F(x)=U^{-1}\bigl(f(x+a_0)-f(a_0)\bigr),\qquad x\in A-a_0,
\]
which is a step-isometry, sends zero to zero, and satisfies
\(\norm{F(x)-x}\le6\).
Proposition~\ref{prop:allbounded} gives \(F(x)=x\). Hence
\[
 f(a)=U(a-a_0)+f(a_0)\qquad(a\in A).
\]
This formula defines a surjective real affine isometry on \(E\).
Any two isometric extensions agreeing on the dense set \(A\) agree
everywhere by continuity.
\end{proof}

\begin{remark}\label{rem:relative}
The smoothness used in the reduction is smoothness of the norm
\emph{restricted to \(Y\)}. It does not assert that \(Y\) contains smooth
points of the ambient space \(E\). Supporting functionals belong to
\(Y^*\), and no smoothness hypothesis on the ambient space is needed.
\end{remark}

\section{Step-isometries on \texorpdfstring{$\ell_p$}{ell1} spaces}\label{sec:lone}
\begin{theorem}\label{main-lone-finite}
Let \(A\subset X\) be a norm-dense subset, and let \(f:A\to X\) be a step-isometry.
\begin{enumerate}
\item[(i)]If $X=\ell_1^n$ with integer $n\ge3$, then $f$ admits the form
\begin{equation}\label{finite-signed}
f(x_1,\ldots,x_n)=(\varepsilon_1x_{\pi(1)},\ldots,\varepsilon_nx_{\pi(n)})+b,
\end{equation}
where $b\in\ell_1^n$, each $\varepsilon_j\in\{-1,1\}$, and $\pi$ is a permutation of \(\{1,2,\dots,n\}\). In particular, property $\ell_1^n$ satisfies Property A.
\item[(ii)] If \(X=\ell_1\), assume additionally that $f$ is almost surjective. Then
    \begin{equation}\label{infinite-signed}
f((x_j)_{j\ge1})=(\varepsilon_jx_{\pi(j)})_{j\ge1}+b,
\end{equation}
where $b\in\ell_1$, $\pi:\N\to\N$ is a permutation, and each $\varepsilon_j\in\{-1,1\}$. Consequently $\ell_1$ satisfies Properties $\mathrm B$ and $\mathrm C$.
\end{enumerate}
\end{theorem}

\begin{proof}
Let \(X=\ell_1(I)\), where \(I=\{1,\ldots,n\}\) with \(n\ge3\) for the finite-dimensional case, and \(I=\mathbb N\) in the
infinite-dimensional case. In the latter case, we further assume that
\(f\) is almost surjective.

Fix \(a_0\in A\). by
Lemma~\ref{lem:normalization}, there exists a surjective real-linear
isometry \(U:X\to X\) and a normalized map
\[
 F(x)=U^{-1}\bigl(f(x+a_0)-f(a_0)\bigr),\qquad x\in A-a_0,
\]
which is a step-isometry, sends 0 to 0, and satisfies
\(\norm{F(x)-x}\le6\).

We aim to prove that
\[
F(x)=x\qquad(x\in A-a_0).
\]
For each \(s=(s_j)_{j\in I}\in\{-1,1\}^I\), define
\[
\phi_s(x)=\sum_{j\in I}s_jx_j.
\]
Then \(\phi_s\in S_{X^*}\).

Choose \(c=\{c_j\}_{j\in I}\in \R^I \) such that $c_j>0$ for all $j\in I$  and
\[
\sum_{j\in I}c_j=1.
\]
Put $u=\{u_j\}_{j\in I}$ where \(u_j=s_jc_j\). Then \(u\in S_X\)  and \(u_j\not=0\) for every $j\in I$. For each \(z=\{z_j\}_{j\in I}\in X\), we have
\[
\lim_{t\to0}
\frac{|u_j+tz_j|-|u_j|}{t}=s_jz_j,
\qquad
\left|
\frac{|u_j+tz_j|-|u_j|}{t}
\right|\le |z_j|.
\]
Since \(\sum_{j\in I}|z_j|<\infty\), we deduce that
\[
\lim_{t\to0}
\frac{\|u+tz\|-\|u\|}{t}
=
\sum_{j\in I}s_jz_j
=
\phi_s(z).
\]
Consequently, the norm is G\^ateaux differentiable at \(u\),
with its G\^ateaux derivative \(j(u)=\phi_s\).
By Proposition \ref{prop:order}, it follows that
\[
\phi_s(Fx)\le\phi_s(Fy)
\]
for all $x,y\in A-x_0$ satisfying $\phi_s(x)<\phi_s(y)$.

Fix \(s\in\{-1,1\}^I\). Since \(|I|\ge3\), we may take a partition
\(
I=I_1\sqcup I_2\sqcup I_3
\)
into three nonempty subsets. Define
\[
r_i(x)=\sum_{j\in I_i}s_jx_j
\qquad(i=1,2,3),
\]
and set
\[
\begin{aligned}
\psi_1&=-r_1+r_2+r_3,\\
\psi_2&=r_1-r_2+r_3,\\
\psi_3&=r_1+r_2-r_3.
\end{aligned}
\]
It is clear that \(\psi_1,\psi_2,\psi_3\) are linearly independent, and
\[
\phi_s=\psi_1+\psi_2+\psi_3.
\]
Since \(\{\phi_s, \psi_1,\psi_2,\psi_3\}\subset \Omega_F\), Lemma \ref{lem:finitecontinuous} yields
\[
\psi_i(Fx)=\psi_i(x),\qquad
\phi_s(Fx)=\phi_s(x)
\qquad i=1,2,3;\ x\in A-x_0.
\]
Recall that
\[
\|z\|=\sup_{s\in\{-1,1\}^I}|\phi_s(z)|
\qquad(z\in X).
\]

Thus for any \(x\in A-x_0\),
\[
\begin{aligned}
\|Fx-x\|
&=\sup_s|\phi_s(Fx-x)|\\
&=\sup_s|\phi_s(Fx)-\phi_s(x)|\\
&=0.
\end{aligned}
\]
Consequently,
\(
F=I|_{A-x_0}.
\)
From the definition of \(F\),
\[
f(a_0+x)-f(a_0)=Ux
\qquad(x\in A-x_0).
\]
Thus, for every \(x\in A\),
\[
f(x)=U(x-a_0)+f(a_0)=Ux+b,
\qquad
b=f(a_0)-Ua_0.
\]
In particular, \(f\) is the restriction of a surjective
affine isometry of \(X\).

\end{proof}

\begin{remark}
 $\ell_1^2$ does not satisfy $\mathrm Property\quad A$, $\mathrm B$, and $\mathrm C$ in [7], so \(X=\ell_1^n\) in Theorem~\ref{main-lone-finite} holds for $n\ge 3$.
\end{remark}

The following example demonstrates that the almost-surjectivity hypothesis in Theorem~\ref{main-lone-finite} cannot be dropped in the case \(X=\ell_1\).

\begin{example}
Define \(H:\ell_1\to\ell_1\) by
\[
H(x)=(x_1^+,x_1^-,x_2^+,x_2^-,\ldots),
\qquad \forall x=(x_n)_{n\in \N}\in \ell_1,
\]
where$$t^+=\max\{t,0\},\quad t^-=\max\{-t,0\},$$ for every $t\in \R$.
The identity
\[
|s^+-t^+|+|s^--t^-|=|s-t|
\]
implies that
\[
\|H(x)-H(y)\|=\|x-y\|.
\]
Hence \(H\) is a step-isometry.  Since every coordinate
of \(H(x)\) is nonnegative, we have
\[
\operatorname{dist}(-te_1,H(\ell_1))=t
\qquad(t>0),
\]
which shows that \(H\) is not almost surjective.

On the other hand,
\[
H(0)=0,\qquad H(e_1)=e_1,\qquad H(-e_1)=e_2,
\]
so \(H\) is not affine.
\end{example}

\begin{remark}
The standard space \(\ell_2\) shows that uniform nonsquareness
does not imply \(\mathrm Property\) A.

Indeed, for any \(a,b\in S_{\ell_2}\), the parallelogram
identity gives
\[
\|a+b\|_2^2+\|a-b\|_2^2=4.
\]
Consequently,
\[
J(\ell_2)
=
\sup_{a,b\in S_{\ell_2}}
\min\{\|a+b\|_2,\|a-b\|_2\}
\le\sqrt2<2,
\]
so \(\ell_2\) is uniformly non-square.

Consider the right shift
\[
S:\ell_2\to\ell_2,\qquad
S(x_1,x_2,\ldots)=(0,x_1,x_2,\ldots).
\]
This map is linear and satisfies
\[
\|Sx-Sy\|_2^2
=
\sum_{j=1}^{\infty}|x_j-y_j|^2
=
\|x-y\|_2^2.
\]
Thus \(S\) is an isometric embedding and hence a
step-isometry. However,
\[
S(\ell_2)=\{x\in\ell_2:x_1=0\},
\]
so \(S\) is not surjective. Taking \(A=\ell_2\), any affine
isometry extending \(S|_A\) must equal \(S\) on the whole
space and therefore cannot be surjective. Hence
\(\mathrm P_1(\ell_2)\) fails.

This does not contradict \(\mathrm P_2\) or \(\mathrm P_3\):
the range of \(S\) is not dense, and
\[
\operatorname{dist}(te_1,S(\ell_2))=|t|,
\]
so \(S\) is not almost surjective.
\end{remark}

The following result from~\cite[Proposition 2.3]{SV}.
\begin{lemma}\label{asymptotic}
Let $E,F$ be real Banach spaces, with $F$ uniformly convex, and let $h:E\to F$ be an $\varepsilon$-isometry with $h(0)=0$. Then
\[
Ux=\lim_{t\to\infty}\frac{h(tx)}t
\]
exists in norm and defines a linear isometric embedding.
\end{lemma}

\begin{lemma}\label{projection}
Let $U:\ell_p\to\ell_p$ be a  linear isometric embedding, where $1<p<\infty$. then there exists a norm-one projection $P$ from $\ell_p$ onto $U(\ell_p)$.
\end{lemma}
\begin{proof}
By \cite[Theorem 2.4 and Remark 2.5(1)]{SV}, there exists a linear operator
\(T:\ell_p\to\ell_p\) such that
\[
\|T\|=1,\qquad TU=I.
\]
Set \(P=UT\). Then
\[
P^2=U(TU)T=P,
\qquad
P(\ell_p)=U(\ell_p),
\]
where the equality for the range follows from the identity \(PU=U\).
Since \(U\) is an isometry,
\[
\|Px\|_p=\|Tx\|_p
\qquad(x\in\ell_p),
\]
which yields \(\|P\|=\|T\|=1\).

\end{proof}

\begin{theorem}\label{lp}
Let $1<p<\infty$, let $A$ be dense subset of $\ell_p$, and let   $f:A\to\ell_p$ be a step-isometry. Then $f$ is the restriction of a unique affine isometric embedding
\[
f(x)=Ux+b\qquad(x\in A).
\]
If $p=2$, $U$ may be any linear isometric embedding. If $p\ne2$, the representation takes the form
\[
f(x)=b+\sum_{i=1}^{\infty}x_i u_i,
\]
where $(u_i)$ is a sequence of disjointly supported unit vectors in $\ell_p$.

Moreover, the almost surjectivity of $f$, the density of $f(A)$, and surjectivity of $U$ are equivalent.
\end{theorem}

\begin{proof}
We may assume that \(0\in A\) and \(f(0)=0\). As in the proof of  Lemma \ref{lem:normalization},  there exists a \(3\)-isometry \(h:\ell_p\to \ell_p\) such that
\[
h|_A=f, h(0)=0.
\]

By Lemmas \ref{asymptotic} and  \ref{projection},
$$
Ux=\lim_{t\to\infty}\frac{h(tx)}{t}
$$
defines a linear isometric embedding \(U:\ell_p\to \ell_p\), and there exists a norm-one projection \(P\) from $\ell_p$ onto \(U(\ell_p)\). Let \(Q=U^{-1}P\), then
\begin{equation}\label{PQ}
\|Q\|=1,\qquad QU=I,\qquad UQ=P.
\end{equation}

Since the norm of \(\ell_p\) is G\^ateaux differentiable
away from the origin and
$$
j(S_{\ell_p})=S_{\ell_p^*}.
$$
For any \(u\in S_{\ell_p}\), we have \(\|j(Uu)P\|\le 1\) and  \(j(Uu)P(Uu)=\|Uu\|\). The uniqueness of the
supporting functional at \(Uu\) therefore yields
$$
j(Uu)P=j(Uu).
$$
Similarly,   \(j(Uu)U=j(u)\). Consequently,
\begin{equation}\label{supportP}
j(Uu)=j(Uu)P=j(Uu)UQ=j(u)Q.
\end{equation}

We now  prove separately that \(Qf=I|_A\) and \((I-P)f=0\).

Fix \(u\in S_{\ell_p}\) and \(x,y\in A\) satisfying
\(
j(u)(x)<j(u)(y).
\)
Choose \(c\in \R\) and \(\{z_n\}_{n\in \N}\subset A\) such that
$$
j(u)(x)<c<j(u)(y),\quad and \quad \|z_n-(n+c)u\|<\frac1n, \quad \forall n\in\N.
$$
Using the same source-side expansion as in Proposition \ref{prop:order}, there exists $N_0\in \mathbb N$ such that
\[
\|fz_n-fx\|\ge n>\|fz_n-fy\|,
\]
for all $n\ge N_0$.
Moreover, since \(h\) extends \(f\) and is a \(3\)-isometry,
\[
\|fz_n-h(nu)\|
\le \|z_n-nu\|+3
\le |c|+\frac1n+3.
\]
Hence,
$$
\lim\limits_{n\to\infty}\frac{fz_n}{n}=\lim\limits_{n\to\infty}\frac{h(nu)}{n}= Uu.
$$
By Lemma \ref{lem:fixed},
\[
j(Uu)(fy-fx)=\lim\limits_{n\to\infty}(\|fz_n-fx\|-\|fz_n-fy\|)\ge0.
\]
Since \(j(S_{\ell_p})=S_{\ell_p^*}\), for any $\phi\in S_{\ell_p^*} $, there exists $u\in S_{\ell_p}$ with $j(u)=\phi$. Combining (\ref{supportP}) and theargument above,   we obtain
\begin{equation}\label{allorder}
\phi(Qfy-Qfx)=j(u)Q(fy-fx)=j(Uu)(fy-fx)\ge0
\end{equation}
for all $x,y\in A$ with $\phi(x-y)>0$.

It follows that for any $x,y\in A$ with $x\ne y$, there exists $\lambda_{xy}\ge0$ such that
\[
Qfx-Qfy=\lambda_{xy}(x-y).
\]
Since \(Qf(0)=0\), there exist $\lambda_{x}, \lambda_{y}\ge0$ satisfying
\[
Qfx=\lambda_xx,
\qquad Qfy=\lambda_yy,
\]
for all \(x,y\in A\setminus\{0\}\).

If \(x,y\in A\) are linearly
independent, then
\[
\lambda_xx-\lambda_yy
=Qfx-Qfy
=\lambda_{xy}(x-y).
\]
Linear independence forces
\[
\lambda_x=\lambda_{xy}=\lambda_y.
\]
If \(x,y\) are nonzero and collinear, since $A$ is dense in $\ell_p$ and
\(\dim \ell_p\ge2\),  there exists \(z\in A\) such that $z\not\in \R x$.
Comparing with \(z\) yields the same conclusion. Therefore
there exists a   constant \(\lambda\ge0\) such that
\begin{equation}\label{xishu}
Qfx=\lambda x,\qquad \forall x\in A.
\end{equation}

Fix \(u\in S_{\ell_p}\), and pick \(w_n\in A\) such that
\[
\|w_n-nu\|<\frac1n.
\]
By the same reasoning as before,
$$
\lim\limits_{n\to\infty}\frac{fw_n}{n}= Uu.
$$
Applying (\ref{PQ}) and (\ref{xishu}),
\[
\lambda u
=\lim_{n\to\infty}\frac{gw_n}{n}
=QUu=u.
\]
Thus \(\lambda=1\),  so
\begin{equation}\label{projection-exact}
Qfx=x,\qquad Pfx=Ux
\qquad(x\in A).
\end{equation}
Fix distinct \(x,y\in A\), and let $k\in \N$ with
\(
k>\frac{\|x-y\|}{2}.
\)
Then there exists \(z_*\in X\) such that
\[
\|z_*-x\|=\|z_*-y\|=k.
\]
Set
$$
v=\frac{y-x}{\|y-x\|}\quad and \quad m=\frac{x+y}2.
$$
Then
\[
\|x+kv-y\|=|k-\|x-y\||<k
\quad\text{and}\quad \|x-kv-y\|=
k+\|x-y\|>k.
\]

For \(0<t<1\),
\[
\|(1-t)z_*+tm-x\|=\|(1-t)(z_*-x)+\frac{t(y-x)}{2}\|
\le (1-t)k+t\frac{\|x-y\|}{2}<k,
\]
and \[
\|(1-t)z_*+tm-y\|=\|(1-t)(z_*-y)+\frac{t(x-y)}{2}\|
\le (1-t)k+t\frac{\|x-y\|}{2}<k.
\]
Since $A$ is dense in $\ell_p$, we may select
\[
\{z_n\}_{n\in \N}\subset A\cap B(x,k)\cap B(y,k),
\]
such that $\lim\limits_{n\to \infty}z_n=z_*.$ Let
\[
a_n=fx-fz_n,\qquad b_n=Pa_n.
\]
Since $\|x-z_n\|<k$,  $f$ is step-isometry, and $P$ is a norm-one projection, it follows that
\begin{equation}\label{anhebnfanshu}
\|a_n\|<k\quad and \quad \|b_n\|<k.
\end{equation}
By (\ref{projection-exact}),
\begin{equation}\label{bnfanshu}
\lim\limits_{n\to\infty}\|b_n\|=\lim\limits_{n\to\infty}\|U(x-z_n)\|=\lim\limits_{n\to\infty}\|x-z_n\|= k.
\end{equation}
Since \(P^2=P\) and \(\|P\|=1\), together with (\ref{bnfanshu}), we get
\begin{equation}\label{hedefanshu}
\lim\limits_{n\to\infty}\left\|\frac{a_n+b_n}{2}\right\|
\ge
\lim\limits_{n\to\infty}\left\|P\left(\frac{a_n+b_n}{2}\right)\right\|
=\lim\limits_{n\to\infty}\|b_n\|= k.
\end{equation}
The space \(l_p\) is   uniform convex, so (\ref{anhebnfanshu}) and (\ref{hedefanshu}) imply
\[
\lim\limits_{n\to\infty}\|a_n-b_n\|=0.
\]
Therefore,
\[
\lim\limits_{n\to\infty}\big((I-P)f(x)-(I-P)f(z_n)\big)=\lim\limits_{n\to\infty}(a_n-b_n)=0.
\]
Applying the same argument to \(y\), using the same sequence
\(z_n\), yields
\[
\lim\limits_{n\to\infty}\big((I-P)f(x)-(I-P)f(z_n)\big)=0.
\]
Hence \((I-P)f(x)=(I-P)f(y)\). That is, \((I-P)f\) is constant on \(A\). Since \((I-P)f(0)=0\), $$(I-P)f(x)=0,\quad \forall x\in A.$$ Combining this with (\ref{projection-exact}), we conclude
\[
f(x)=Pf(x)=Ux,\qquad \forall x\in A.
\]
Uniqueness of the affine isometric extension follows from density.

Finally, since \(A\) is dense and \(U\) is an isometry,
\[
\overline{f(A)}=U(\ell_p).
\]
It follows that density of \(f(A)\) is equivalent to surjectivity of
\(U\). If \(U(\ell_p)\ne \ell_p\), choose \(z\notin U(\ell_p)\).
The range \(U(\ell_p)\) is closed, so
\[
\lim\limits_{t\to\infty}\operatorname{dist}(tz,f(A))
=\lim\limits_{t\to\infty}t\,\operatorname{dist}(z,U(X))
=\infty.
\]
Thus almost surjectivity   forces \(U\) to be surjective.
Conversely, if \(U\) is surjective, then \(f(A)\) is dense and
hence almost surjective.

This completes the proof.

\end{proof}

\section{The geometric rigidity theorems}\label{sec:geometry}
In this section, we discuss the rigidity properties for general Banach spaces, addressing Questions \ref{question4}, \ref{question5} and \ref{question6} posed by \v{S}emrl \cite{Semrl}.

Let $X$ be a real Banach space, and let $A$ be a dense subset of $X$.
A step-isometry \(F:A\to X\) is said to be \textbf{normal}  if it satisfies
\[
  F(0)=0,\qquad \sup_{x\in A}\norm{F(x)-x}<\infty.
\]

\begin{lemma}\label{lianggefanhan}
Let $X$ be a real Banach space, and let $A$ be a dense subset of $X$, and suppose
\(F:A\to X\) is a normal  step-isometry. If $ \{\phi,\psi\}\in \Omega_F$ and $\phi,\psi$ are linearly independent, then
\begin{equation}\label{functional-estimate}
| \phi(Fx-x)|\le M\norm{ \phi-\psi},\qquad \forall x\in A.
\end{equation}
\end{lemma}
\begin{proof}
Fix \(x\in A\) and \(\epsilon >0\), and set \(t = \phi (x)\). Since \(\phi\) and \(\psi\) are linearly independent, the  mapping
$$\Phi(z)= (\phi (z),\psi (z)),\qquad \forall z\in X$$
 is surjective onto \(\mathbb{R}^{2}\). Let $$R_+ = (t, t+\epsilon) \times (-\epsilon, 0), \quad R_- = (t-\epsilon, t) \times (0, \epsilon).$$ Both \(R_+\) and \(R_-\)  are nonempty open subsets of \(\mathbb{R}^{2}\). By   density of \(A\), we may select \(z_{+},z_{- }\in A\) such that
\[
t < \phi(z_+) < t+\epsilon, \quad -\epsilon < \psi(z_+) < 0,
\]
\[
t-\epsilon < \phi(z_-) < t, \quad 0 < \psi(z_-) < \epsilon.
\]
Let \(d(z) = Fz - z\). Since \(F(0) = 0\) and $ \{\phi,\psi\}\in \Omega_F$,
\[
\phi (Fz_{-})\leq \phi (Fx)\leq \phi (Fz_{+}), \qquad \psi (Fz_{+})\leq 0\leq \psi (Fz_{-}).
\]
We estimate \(\phi(Fx) - t\) from above:
\begin{align*}
\phi (Fx) - t
&\le \phi (Fz_{+}) - t \\
&= \phi (z_{+}) - t + \phi(d(z_{+})) \\
&= (\phi (z_{+}) - t) + \psi(d(z_{+})) + (\phi -\psi)d(z_{+}) \\
&\le (\phi (z_{+}) - t) - \psi(z_{+}) + (\phi -\psi)d(z_{+}) \\
&< 2\epsilon + M\|\phi -\psi \|.
\end{align*}
Similarly,
\begin{align*}
\phi (Fx) - t
&\ge \phi (Fz_{-}) - t \\
&= \phi (z_{-}) - t + \phi(d(z_{-})) \\
&= (\phi (z_{-}) - t) + \psi(d(z_{-})) + (\phi -\psi)d(z_{-}) \\
&\ge (\phi (z_{-}) - t) - \psi(z_{-}) + (\phi -\psi)d(z_{-}) \\
&> -2\epsilon -M\|\phi -\psi \|.
\end{align*}
Hence $$|\phi(Fx) - t| < 2\epsilon + M\|\phi - \psi\|.$$ Since \(\epsilon > 0\) is arbitrary, we obtain
$$|\phi(Fx) - t| \le  M\|\phi - \psi\|.$$

This completes the proof.

\end{proof}

\begin{proposition}\label{choumi}
Let $X$ be a real Banach space, and let $A$ be a dense subset of $X$, and suppose
\(F:A\to X\) is a normal  step-isometry. Let $\Gamma\subset S_{X^*}$ be a norming family, i.e.,
\[
\norm z=\sup_{\varphi\in\Gamma}|\varphi(z)|\qquad(z\in X).
\]
If $\Gamma\subset\Omega_F$ and
$$\varphi\in \overline{\Gamma\setminus\{\varphi,-\varphi\}}^{\|\cdot\|},$$ for every $\varphi\in\Gamma$, then $F=I|_A$.
\end{proposition}
\begin{proof}
Fix an arbitrary \(\phi \in \Gamma\). By hypothesis, \(\phi\) lies in the norm closure of \(\Gamma \setminus \{\phi, -\phi\}\). Therefore, there exists a sequence \(\{\psi_n\} \subset \Gamma \setminus \{\phi, -\phi\}\) such that $$ \lim\limits_{n \to \infty}\|\psi_n - \phi\| = 0.$$

Since \(\phi\) and \(\psi_n\) are linearly independent for each \(n\).  Lemma \ref{lianggefanhan} yields that for any \(x \in A\),
\[
|\phi (Fx - x)|\leq M\|\phi -\psi_n\|.
\]
Taking \(n \to \infty\) and using \(\|\psi_n - \phi\| \to 0\), we obtain
\[
 \phi(Fx - x) = 0.
\]
 Finally, by the norming property of \(\Gamma\),
\[
\|Fx - x\| = \sup_{\phi \in \Gamma}|\phi (Fx - x)| = 0.
\]
Hence \(F(x) = x\) for all \(x \in A\).

\end{proof}
    In a reflexive Banach space, the G\^ateaux smooth unit vectors of any given norm are dense in $S_X$. This follows, for example, from Asplund's differentiability theorem~\cite{Asplund}; reflexive spaces are Asplund, so one even obtains a dense set of Fr\'echet smooth points.
\begin{lemma}\label{zifanyinli}
Let $X$ be a real reflexive Banach space, and let $A$ be a dense subset of $X$, and suppose
\(F:A\to X\) is a normal  step-isometry. If $F\ne I|_A$, then there exist $v\in S_X$ and $\varphi\in S_{X^*}$ with $\varphi(v)=1$ such that
\begin{equation}\label{flat}
\norm{v+y}=\norm{v-y}=1\qquad( \forall y\in\ker\varphi,\ \norm y\le1).
\end{equation}
Moreover, the same $(v,\varphi)$ satisfies
\begin{equation}\label{max-identity}
\norm x=\max\{|\varphi(x)|,\norm{x-\varphi(x)v}\}\qquad( \forall x\in X).
\end{equation}
\end{lemma}
\begin{proof}
Since \(X\) is reflexive, the G\^ateaux smooth points of $X$ are dense in
\(S_X\). Set
$$\Gamma=\{j(u): u \in s.m.(X)\}.$$ Then \(\Gamma\) is symmetric and norming, i.e.,
\[
\|z\|=\sup_{\psi\in \Gamma}|\psi(z)|
\qquad(z\in X).
\]
Indeed, for \(z\ne0\), choose smooth unit vectors \(u_n\) such that
\(\lim\limits_{n\to \infty}u_n= z/\|z\|\). Then
\[
1\ge \lim\limits_{n\to \infty}j(u_n)\left(\frac{z}{\|z\|}\right)
\ge
1-\lim\limits_{n\to \infty}\left\|u_n-\frac{z}{\|z\|}\right\|
=1.
\]

If \(F\ne I|_A\), there exists \(x_0\in A\) with \(F(x_0)-x_0\ne0\).
By the norming property of \(\Gamma\), we may pick  \(\varphi\in\Gamma\) satisfying
\[
\varphi((F(x_0)-x_0)\ne0.
\]

By Proposition \ref{prop:order}, \(\Gamma\subset \Omega_F\).
We claim that
\begin{equation}\label{notspan}
\varphi\notin \operatorname{span}\bigl(\Gamma\setminus\{\varphi,-\varphi\}\bigr).
\end{equation}
Suppose, for contradiction, that the contrary holds. Then we can write
\(
\varphi=\sum_{i=1}^{m}a_i\psi_i,
\)
where \(\psi_1,\ldots,\psi_m\) are linearly independent elements of
\(\Gamma\setminus\{\pm\varphi\}\)  and every \(a_i\ne0\). Lemma~\ref{lem:finitecontinuous} then implies
\(
\varphi(Fx_0)=\varphi(x_0),
\)
which   contradicts our choice of \(x_0\) and $\varphi$.

The function \(g_\varphi(t)-t\) cannot be constant on the continuity
set of \(g_\varphi\).
Suppose otherwise that its value equals \(c\) there. Then
\[
g_\varphi(t)=t+c\qquad(t\in\mathbb R).
\]
By Lemma~\ref{lem:continuous},
\[
\varphi(Fx)=g_\varphi(\varphi(x))=\varphi(x)+c\qquad(x\in A).
\]
Since \(F(0)=0\), we would obtain \(c=0\), contradicting
\(\varphi((F(x_0)-x_0)\ne0\).

Thus there exist continuity points \(r_1,r_2\) such that
\[
a_k=g_\varphi(r_k)-r_k \quad(k=1,2),
\qquad
\delta:=|a_1-a_2|>0.
\]

We next prove that
\begin{equation}\label{positive-distance}
\delta\le2M\norm{\psi-\chi},\qquad
\left(\quad \forall
\chi\in
\operatorname{span}\bigl(\Gamma\setminus\{\pm\varphi\}\bigr)
\right).
\end{equation}
Take a nonzero such \(\chi\), and write
\(
\chi=\sum_{i=1}^{m}c_i\psi_i
\)
where \(\psi_1,\ldots,\psi_m\) are linearly independent members of
\(\Gamma\setminus\{\pm\varphi\}\).
By (\ref{notspan}),
\(\psi,\psi_1,\ldots,\psi_m\) are linearly independent.
Hence the linear map
\[
\Phi(z)=
\bigl(\psi(z),\psi_1(z),\ldots,\psi_m(z)\bigr), \quad \forall z\in X
\]
is surjective onto \(\mathbb R^{m+1}\).

Choose a continuity point \(s_i\) of each \(g_{\psi_i}\).
Surjectivity of $\Phi$ guarantees the existence of  \(z_1,z_2\in X\) satisfying
\[
\psi(z_k)=r_k,\qquad
\psi_i(z_k)=s_i
\quad(k=1,2,\ 1\le i\le m).
\]
Select \(x_{k,n}\in A\) with \(\lim\limits_{n\to \infty}x_{k,n}= z_k\). By Lemma~\ref{lem:continuous},
\begin{align*}
\lim\limits_{n\to \infty}\psi(Fx_{k,n}-x_{k,n})&= a_k,\\
\lim\limits_{n\to \infty}\psi_i(Fx_{k,n}-x_{k,n})&= g_{\psi_i}(s_i)-s_i.
\end{align*}
On the other hand,   Lemma \ref{lianggefanhan} gives
\[
|(\psi-\chi)(Fx_{k,n}-x_{k,n})|\le M\norm{\psi-\chi}.
\]
Therefore
$$\delta=\lim\limits_{n\to \infty}|(\psi-\chi)(Fx_{1,n}-x_{1,n})-(\psi-\chi)(Fx_{2,n}-x_{2,n})|\le2M\norm{\psi-\chi}.$$
Set
\[
V=
\overline{
\operatorname{span}\bigl(\Gamma\setminus\{\pm\varphi\}\bigr)
}^{\|\cdot\|_{X^*}}.
\]
It follows that
\[
d_0:=\operatorname{dist}(\varphi,V)
\ge\frac{\delta}{2M}>0.
\]
Define a linear functional on \(V+\mathbb R\varphi\) by
\[
L(w+t\varphi)=t.
\]
This definition is well defined, and the inequality
\[
\|w+t\varphi\|\ge |t|d_0
\]
shows that \(L\) is bounded.
By the Hahn--Banach theorem, \(L\) extends to an element
\(\widetilde L\in X^{**}\).
Since \(X\) is reflexive, there exists \(v\in X\) such that
\[
\widetilde L(\psi)=\psi(v)\qquad(\psi\in X^*).
\]
Accordingly,
\[
\varphi(v)=1,\qquad
\psi(v)=0
\quad\bigl(\psi\in\Gamma\setminus\{\pm\varphi\}\bigr).
\]
We obtain
\[
\|v\|
=
\sup_{\psi\in\Gamma}|\psi(v)|
=1.
\]

Now let \(y\in\ker\varphi\).
Since \(\varphi(y)=0\), the norming property yields
\[
\|y\|
=
\sup_{\psi\in\Gamma\setminus\{\pm\varphi\}}|\psi(y)|.
\]
Thus  for every \(t\in\mathbb R\),
\[
\begin{aligned}
\|tv+y\|
&=
\sup_{\psi\in\Gamma}|\psi(tv+y)|\\
&=
\max\left\{
|t|,
\sup_{\psi\in\Gamma\setminus\{\pm\varphi\}}|\psi(y)|
\right\}\\
&=
\max\{|t|,\|y\|\}.
\end{aligned}
\]
In particular, if \(\|y\|\le1\), taking \(t=1\) and replacing
\(y\) by \(-y\) to obtain
\[
\|v+y\|=\|v-y\|=1,
\]
which establishes (\ref{flat}).

Finally, for arbitrary \(x\in X\), set
\[
t=\varphi(x),\qquad y=x-\varphi(x)v.
\]
Since \(\varphi(v)=1\), we have \(y\in\ker\varphi\).
The preceding identity then gives
\[
\|x\|
=
\max\bigl\{
|\varphi(x)|,\,
\|x-\varphi(x)v\|
\bigr\},
\]
which is exactly (\ref{max-identity}).
\end{proof}

\begin{theorem}\label{yizhifeifang}
Suppose that $X$ is uniformly non-square.
Then $X$ has   Property $B$.
\end{theorem}

\begin{proof}
Let $A$ be a dense subset of $X$, and let $f:A\to X$ be almost surjective. By normalizing via  Lemma \ref{lem:normalization}, it suffices to show that $F= I|_A$. Assume, to the contrary, that this is not the case. By James's theorem~\cite{James}, $X$ is reflexive. By Lemma~\ref{zifanyinli}, there exists a pair $(v,\varphi)$ satisfying \eqref{flat}. Choose $y\in\ker\varphi$ with $\norm y=1$, this is possible because $\dim X\ge2$. Set $a=v+y$ and $b=v-y$. Then $a,b\in S_X$, and we have
\[
\norm{a+b}=2\norm v=2,\qquad \norm{a-b}=2\norm y=2,
\]
which contradicts the assumption that $X$ is uniformly non-square.

\end{proof}

\begin{theorem}\label{yizhifeifang}
Suppose that $X$ is a smooth normed space.
Then $X$ has   Property $B$.
\end{theorem}

\begin{proof}
Let $A$ be a dense subset of $X$, and let $f:A\to X$ be almost surjective. Set \[\Gamma = \{j(u) : u \in S_X\}.\] By the smoothness of $X$ and Proposition \ref{prop:order}, \(\Gamma\subset \Omega_F\), and \(\Gamma\) consists of all norm-attaining unit functionals. The Bishop-Phelps theorem implies that \(\Gamma\) is norm-dense in \(S_{X^*}\) and hence norming.

Given $\varphi\in \Gamma$, choose $\eta$ independent of $\varphi$ and set $$\chi_n=(\varphi+\eta/n)/\norm{\varphi+\eta/n}.$$ Then $\lim\limits_{n\to \infty}\chi_n=\varphi$ and $\chi_n\notin\{\pm\varphi\}$. By the density of $\Gamma$, we may choose $\psi_n\in\Gamma$ such that
\[
\norm{\psi_n-\chi_n}<\min\{1/n,\tfrac12 \dist(\chi_n,\{\varphi,-\varphi\})\}.
\]
It follows that $\psi_n\ne\pm\varphi$ and $\lim\limits_{n\to \infty}\psi_n=\varphi$. Proposition \ref{choumi} yields $F=I|_A$, where $F$ is the normalization of $f$ obtained via  Lemma \ref{lem:normalization}. Therefore,  $X$ has   Property $B$.

\end{proof}

The following theorem, which gives a characterization of Fr\'echet differentiable points, is due to V. L. \v{S}mulian. The proof can be found in \cite{Fabian}.

\begin{theorem}[\v{S}mulian Theorem]\label{frechet-support}
Let \(u\in S_X\). The norm is Fr\'echet differentiable at $u$ if and only if the following holds: whenever \(\varphi_n\in B_{X^*}\) and \(\lim\limits_{n\to\infty}\varphi_n(u)= 1\), we have \(\lim\limits_{n\to\infty}\|\varphi_n - \varphi\|=0\), where $\varphi$ is the unique support functional at $u$.
\end{theorem}

\begin{theorem}\label{yangetuyantuo}
Suppose that $X$ is strictly convex and that the Fr\'echet smooth points of its norm are dense in $S_X$.
Then $X$ has   Property $B$.
\end{theorem}
\begin{proof}
Let $A$ be a dense subset of $X$, and let $f:A\to X$ be almost surjective. Let $D\subset S_X$ denote the set of Fr\'echet smooth points, and set $\Gamma=\{j(u):u\in D\}$. Then $\Gamma\in \Omega_F$, and by the density of \(D\), \(\Gamma\) forms a norming family.

Fix \(u \in D\).  There exists a sequence \(\{u_n\} \subset D \setminus \{u\}\) converging to \(u\). Put \(\psi_n = j(u_n)\) and \(\varphi = j(u)\). We have
\[
1\ge \lim\limits_{n\to \infty}\psi_n(u) = \lim\limits_{n\to \infty} (\psi_n(u_n) + \psi_n(u-u_n)) \ge \lim\limits_{n\to \infty}(1 - \|u_n - u\|) = 1.
\]
Theorem \ref{frechet-support} then yields \(\lim\limits_{n\to \infty}\|\psi_n - \varphi\| = 0\).

Strict convexity ensures that \(\psi_n = \varphi\). Thus,  Proposition \ref{choumi} yields $F=I|_A$, where $F$ is the normalization of $f$ obtained via  Lemma \ref{lem:normalization}. Therefore,  $X$ has   Property $B$.
\end{proof}
The following theorem, known as Asplund's differentiability theorem, can be found in
{\cite{Asplund}}
\begin{lemma}\label{lem:asplund-smoothness}
Let $X$ be a real Asplund space, and let $\Omega\subset X$
be nonempty, open and convex. Every continuous convex
function $g:\Omega\to\mathbb R$ is Fr\'echet differentiable
on a norm-dense $G_\delta$ subset of $\Omega$.
\end{lemma}

By Lemma~\ref{lem:asplund-smoothness}, the Fr\'echet smooth points  are dense in $S_X$. Combined with Theorem \ref{yangetuyantuo}, we readily obtain the following theorem.

\begin{theorem}\label{yangetuyantuo1}
\label{extension}
Let $X$ be a strictly convex Banach space with separable dual. For every dense subset $S\subseteq X$
and every bijective step-isometry $T:S\to S$, there
exist a unique surjective linear isometry $U:X\to X$
and a unique vector $b\in X$ such that
\[
T(s)=Us+b\qquad s\in S.
\]
The map
\[
\widehat T:X\to X,
\qquad
\widehat T(x)=Ux+b,\qquad x\in X.
\]
is the unique whole-space bijective step-isometric
extension of $T$.
\end{theorem}

This settles Conjecture~\ref{conj2}, yielding an affine isometric conclusion and removing the countability assumption on the dense set.

\begin{theorem}
Suppose that $X$ is infinite-dimensional, reflexive, and has UKK.
Then $X$ has   Property $B$.
\end{theorem}

\begin{proof}
Let $A$ be a dense subset of $X$, and let $f:A\to X$ be almost surjective. By normalizing via  Lemma \ref{lem:normalization}, it suffices to show that $F= I|_A$. If $F\ne I|_A$, Lemma \ref{zifanyinli} yields a pair $(v,\varphi)$  satisfing \eqref{flat}. The subspace $\ker\varphi$ is closed and infinite-dimensional. Since $X$ is   reflexive, By Riesz's theorem and Banach-Alaoglu Theorem, we may select a sequence $y_n\in S_X\cap\ker\varphi$ such that \(y_n\rightharpoonup 0\) weakly and $\|y_n-y_m\|\ge1/2$ whenever $n\not= m$. Then $v+y_n$ converges to $v$ weakly and
\[
v+y_n\in S_X,\qquad 
\norm{(v+y_n)-(v+y_m)}>1/2.
\]
for all $n,m\in \N$ with $n\not= m$.
Since $X$ has the UKK, there exists  $\delta>0$ such that $\norm v\le1-\delta$, which  contradicts $\norm v=1$.

\end{proof}
The infinite-dimensional assumption is necessary here. In finite dimensions, every weakly convergent sequence converges in norm, so UKK imposes no constraint; the maximum norm on \(\mathbb R^2\) still admits nonisometric bijective step-isometries.

\begin{lemma}\label{hanshu}
Let $\tau:\R\to \R$ be defined by
\[
\tau(t)=\floor t+(t-\floor t)^2\qquad(t\in\R).
\]
Then $\tau$ is an increasing homeomorphism satisfying $\tau(t+k)=\tau(t)+k$ for every $k\in\mathbb Z$, and $\tau$ preserves integer parts of distances. Moreover,  $\tau(0)=0$ and $\tau(1/2)=1/4$.
\end{lemma}
\begin{proof}
Clearly, $\tau$ maps each interval $[n,n+1]$ increasingly and continuously onto itself. The union over these intervals yields bijectivity, continuity of the inverse, and the translation identity.

For any $s,t\in \R$,  assume without loss of generality that  $s<t$, and set $m=\floor{t-s}$. Then $s+m\le t<s+m+1$. Monotonicity together with the translation identity gives
\[
\tau(s)+m\le\tau(t)<\tau(s)+m+1.
\]
Thus $\floor{\tau(t)-\tau(s)}=m$. 
\end{proof}

\begin{theorem}\label{main-reflexive}
Let $X$ be reflexive, the following assertions are equivalent:
\begin{enumerate}
\item $X$ has Property B;
\item $X$ has Property C;
\item every bijective step-isometry $T:S\to S$, for every dense $S\subset X$, is an isometry;
\item every bijective step-isometry $T:X\to X$ is an isometry;
\item there do not exist $v\in S_X$ and $\varphi\in S_{X^*}$, $\varphi(v)=1$, such that
\begin{equation}\label{max-identity}
\norm x=\max\{|\varphi(x)|,\norm{x-\varphi(x)v}\}\qquad(x\in X).
\end{equation}
\end{enumerate}
\end{theorem}

\begin{proof}
It is clear  that 
\(
\textup{(1)}
\Longrightarrow
\textup{(2)}
\Longrightarrow
\textup{(3)}
\Longrightarrow
\textup{(4)}
.
\)
The implication $\textup{(5)}
\Longrightarrow
\textup{(1)}$ follows from Lemma \ref{zifanyinli}.

It remains to prove \textup{(4)}\(\Rightarrow\)\textup{(5)}
by contrapositive.
Suppose that \textup{(5)} fails. Then there exist
\(v\in S_X\) and \(\varphi\in S_{X^*}\), with
\(\varphi(v)=1\) such that
\[
\|x\|
=
\max\bigl\{
|\varphi(x)|,\,
\|x-\varphi(x)v\|
\bigr\}
\qquad(x\in X).
\]
Define
\[
Qx=x-\varphi(x)v.
\]
Then
\[
Qv=0,\qquad
\varphi(Qx)=0,\qquad
Q^2=Q,\qquad and \quad
\|x\|=\max\{|\varphi(x)|,\|Qx\|\}.
\]
Define
\[
T(x)=Qx+\tau(\varphi(x))v
\qquad(x\in X).
\]
where \(\tau\) is the function \(\tau\) from Lemma \ref{hanshu}.
The identities for \(Q\) yield
\[
\varphi(Tx)=\tau(\varphi(x)),
\qquad
Q(Tx)=Qx.
\]
It follows that \(T\) is bijective, with inverse
\[
T^{-1}(y)=Qy+\tau^{-1}(\varphi(y))v.
\]
Indeed, substitution into either composition yields
\(Qx+\varphi(x)v=x\).

For arbitrary \(x,y\in X\), we obtain
\[
\begin{aligned}
\lfloor\|Tx-Ty\|\rfloor
&=
\max\Bigl\{
\lfloor|\tau(\varphi(x))-\tau(\varphi(y))|\rfloor,\,
\lfloor\|Qx-Qy\|\rfloor
\Bigr\}\\
&=
\max\Bigl\{
\lfloor|\varphi(x)-\varphi(y)|\rfloor,\,
\lfloor\|Qx-Qy\|\rfloor
\Bigr\}\\
&=
\lfloor\|x-y\|\rfloor.
\end{aligned}
\]
since
\(
\lfloor\max\{a,b\}\rfloor
=
\max\{\lfloor a\rfloor,\lfloor b\rfloor\},
(a,b\ge0).
\)
Thus \(T\) is a bijective step-isometry on the whole-space. By Lemma \ref{hanshu}, we have
\(
T(0)=0,
T(\frac v2)=\tau(\frac12)v=\frac v4,
\)
so
\[
\|T(\frac v2)-T(0)\|=\frac14
\ne\frac12=\|\frac v2\|.
\]
Therefore \(T\) is not an isometry, and \textup{(4)} fails.
This proves \textup{(4)}\(\Rightarrow\)\textup{(5)}.

\end{proof}

\begin{theorem}\label{main-pone}
 Let $X$ be a uniformly convex Banach space. Then $X$ has the $\mathrm Property$ A  if and only if every linear isometric embedding $U:X\to X$ is surjective.
\end{theorem}
\begin{proof}
Suppose that $X$ has the $\mathrm Property$ A. Let \(U:X\to X\) be a linear isometric embedding. Taking \(A=X\)
and \(f=U\), we obtain a step-isometry on the whole space. By the \(\mathrm P_1(X)\), the map \(U\) itself must be a
surjective affine isometry. Hence \(U\) is surjective.

Conversely, assume that every linear isometric embedding
of \(X\) into itself is surjective. Let \(A\subset X\) be
dense, and let \(f:A\to X\) be a step-isometry.

Choose \(a_0\in A\), and define
\[
A_0=A-a_0,\qquad
f_0(x)=f(a_0+x)-f(a_0)
\quad(x\in A_0).
\]
Then \(A_0\) is dense, \(0\in A_0\), \(f_0(0)=0\), and
\(f_0\) remains a step-isometry. By the argument of   Lemma \ref{lem:normalization}, there exists
a \(3\)-isometry \(h:X\to X\) such that
\[
h(0)=0,\qquad h|_{A_0}=f_0,\qquad h(X)=f_0(A_0).
\]
Since \(X\) is uniformly convex, Lemma \ref{asymptotic} yields a linear isometric embedding \(U:X\to X\) given by
\[
Ux=\lim_{t\to\infty}\frac{h(tx)}{t}.
\]
By hypothesis, \(U\) is surjective. Define
\(
H=U^{-1}h:X\to X.
\)
Then \(H\) is a \(3\)-isometry, \(H(0)=0\), and
\[
\lim\limits_{t\to \infty}\frac{H(tu)}{t}= u.
\]

We claim that
\[
\|H(x)-x\|\le6
\qquad(x\in X).
\]
Fix \(x\in X\) and a G\^ateaux smooth point \(u\in S_X\). We have
\[
\left|
\|H(nu)-H(x)\|-\|nu-x\|
\right|\le3
\]
and
\[
\left|\|H(nu)\|-n\right|\le3.
\]
Consequently,
\[
\left|
\bigl(\|H(nu)-H(x)\|-\|H(nu)\|\bigr)
-
\bigl(\|nu-x\|-n\bigr)
\right|\le6.
\]
Since \(H(nu)/n\to u\), Lemma \ref{lem:fixed} gives
\[
\lim\limits_{n\to\infty}\|H(nu)-H(x)\|-\|H(nu)\|
= -j(u)(H(x)).
\]
Applying the same lemma to \(nu\), we obtain
\[
\lim\limits_{n\to\infty}\|nu-x\|-n= -j(u)(x).
\]
 We conclude that
\[
|j(u)(H(x)-x)|\le6.
\]

Uniform convexity  implies reflexivity, so the G\^ateaux smooth unit vectors are dense in \(S_X\), and
their supporting functionals form a norming family. This yields
\[
\|H(x)-x\|\le6
\qquad(x\in X),
\]
Since \(h=UH\),
it follows that
\[
\|h(x)-Ux\|\le6
\qquad(x\in X).
\]
 The map \(U\) is surjective, hence for each \(y\in X\), there exist \(x\in X\) such that  \(y=Ux\). Therefore
\[
\operatorname{dist}(y,f_0(A_0))
\le\|y-h(x)\|
=\|Ux-h(x)\|
\le6.
\]
Hence
\[
\sup_{y\in X}\operatorname{dist}(y,f_0(A_0))\le6.
\]
We know that \(f_0\) is almost surjective.

A uniformly convex Banach space is uniformly non-square. By Theorem~\ref{yizhifeifang}, there exist a surjective linear isometry \(V:X\to X\)
and a vector \(c\in X\) such that
\[
f_0(x)=Vx+c
\qquad(x\in A_0).
\]
Since \(f_0(0)=0\), we have \(c=0\). We conclude that
\[
f(x)
=V(x-a_0)+f(a_0)
=Vx+b
\qquad(x\in A),
\]
where \(b=f(a_0)-Va_0\). Thus \(f\) is the restriction
of a surjective affine isometry of \(X\), which proves
$\mathrm Property$ A.
\end{proof}

\section{Continuity of whole-space bijective step-isometry}\label{sec:continuity}
We first establish the closedness and convexity needed for a Baire argument.
The subsequent scalar factorization takes place on an invariant countable
set inside its own separable closed ambient subspace.

\begin{lemma}\label{closedballs}
Let $X$ be a real normed space, and let $T:X\to X$ be a bijective step-isometry. Then
\[
T\bigl(\overline B(x,k)\bigr)
=
\overline B(Tx,k)
\qquad x\in X,\ k\in\mathbb N.
\]
Moreover, $T^{-1}(\overline B(y,r))$ is closed and
convex for every $y\in X$ and $r\ge0$.
\end{lemma}
\begin{proof}
 Since $T$ is a bijective step-isometry, we have
\[
T(B(x,k))=B(Tx,k)
\qquad x\in X,\ k\in\mathbb N.
\]
For any $k\in\mathbb N$, we first show that $\|x-y\|\le k$
if and only if 
$B(y,1)\subseteq B(x,k+1).$

If $\|x-y\|\le k$,  it is immediate that \(B(y,1)\subseteq B(x,k+1)\).

Conversely, suppose that $d:=\|x-y\|>k$.
Choose $t$ satisfying
\[
\max\{0,k+1-d\}<t<1
\]
and set
\[
z=y+t\frac{y-x}{d}.
\]
Then
\[
\|z-y\|=t<1,
\qquad
\|z-x\|=d+t>k+1.
\]
Thus $z\in B(y,1)$ but $z\notin B(x,k+1)$, which yields a contradiction.

 It follows that
\[
\begin{aligned}
\|x-y\|\le k
&\Longleftrightarrow
B(y,1)\subseteq B(x,k+1)\\
&\Longleftrightarrow
T(B(y,1))\subseteq T(B(x,k+1))\\
&\Longleftrightarrow
B(Ty,1)\subseteq B(Tx,k+1)\\
&\Longleftrightarrow
\|Tx-Ty\|\le k.
\end{aligned}
\]
Therefore,
\[
T\bigl(\overline B(x,k)\bigr)
=
\overline B(Tx,k).
\]
In particular,
\begin{equation}\label{integer-ball-preimage}
T^{-1}\bigl(\overline B(z,k)\bigr)
=
\overline B(T^{-1}z,k)
\qquad(z\in X,\ k\in\mathbb N).
\end{equation}

Now fix $y\in X$ and $r\ge0$, and select an integer
$m>r$. We claim that
\begin{equation}\label{ballintersection}
\overline B(y,r)
=
\bigcap_{z\in\overline B(y,m-r)}
\overline B(z,m).
\end{equation}
Indeed, if $w\in\overline B(y,r)$ and
$z\in\overline B(y,m-r)$, then
\[
\|w-z\|
\le \|w-y\|+\|y-z\|
\le r+(m-r)=m.
\]
This establishes one inclusion.

Conversely, suppose that
$d:=\|w-y\|>r$, and let
\(
z=y-(m-r)\frac{w-y}{d}.
\)
Then
\[
\|z-y\|=m-r,
\qquad
\|w-z\|=d+m-r>m.
\]
Thus $z\in\overline B(y,m-r)$ but
$w\notin\overline B(z,m)$, which proves
\eqref{ballintersection}.

Taking preimages under $T$ in identity~\eqref{ballintersection}
and applying~\eqref{integer-ball-preimage}, we obtain
\[
\begin{aligned}
T^{-1}\bigl(\overline B(y,r)\bigr)
&=
\bigcap_{z\in\overline B(y,m-r)}
T^{-1}\bigl(\overline B(z,m)\bigr)\\
&=
\bigcap_{z\in\overline B(y,m-r)}
\overline B(T^{-1}z,m).
\end{aligned}
\]
The right-hand side is an intersection of closed
convex sets, and it is itself closed and convex.

\end{proof}
\begin{lemma}\label{barrel}
Let $X$ be a Banach space, and let $F:X\to X$
satisfy $F(0)=0$. Suppose that
\(
F^{-1}\bigl(\overline B(0,r)\bigr)
\)
is closed and convex for every $r>0$. If
\[
\lim_{n\to\infty}\|F(\frac vn)\|=0
\qquad(v\in X),
\]
then $F$ is continuous at 0.
\end{lemma}
\begin{proof}
For any $\varepsilon>0$, choose $0<r<\varepsilon$ and 
define
\[
C_r
=
\{x\in X:\|F(x)\|\le r,\ \|F(-x)\|\le r\}.
\]
Since \(
F^{-1}\bigl(\overline B(0,r)\bigr)
\)
is closed and convex, the set $C_r$ is closed, convex and symmetric.
Moreover, $0\in C_r$ because $F(0)=0$.

We claim that
\begin{equation}\label{absorbing-cover}
X=\bigcup_{n=1}^{\infty}nC_r,
\qquad
nC_r:=\{nx:x\in C_r\}.
\end{equation}
Indeed, for each $v\in X$, applying the hypothesis 
to   $v$ and $-v$ yields
\[
\lim_{n\to\infty}\|F(\frac vn)\|=0,
\qquad
\lim_{n\to\infty}\|F(-\frac vn)\|=0.
\]
Therefore, for all sufficiently large $n$,
\[
\frac vn\in C_r,
\qquad\text{and hence}\qquad
v\in nC_r.
\]
This establishes~\eqref{absorbing-cover}.

Since $X$ is a Banach spaces,  the Baire category theorem implies $C_r$ has nonempty interior. Thus there exist $a\in X$ and $\delta>0$
such that
\[
B(a,\delta)\subseteq C_r.\qquad
B(-a,\delta)\subseteq C_r.
\]

We next verify that $B(0,\delta)\subseteq C_r$.
If $\|h\|<\delta$, then
\[
a+h\in B(a,\delta)\subseteq C_r,
\qquad
-a+h\in B(-a,\delta)\subseteq C_r.
\]
By convexity,
\[
h=\frac{(a+h)+(-a+h)}2\in C_r.
\]
It follows that for any $x\in X$ with $\|x\|<\delta$,  
\[
\|F(x)-F(0)\|
=\|F(x)\|
\le r
<\varepsilon.
\]
Consequently, $F$ is continuous at 0.
\end{proof}

We recall Mazur's differentiability theorem in the form required for our purposes.
\begin{lemma}[{Mazur, \cite{Mazur}}]\label{mazur-smoothness}
For a separable real Banach space $Z$, the G\^ateaux smooth points of the
given norm are norm-dense in $S_Z$.
\end{lemma}

\begin{lemma}\label{invariant-model}
Let $X$ be a real Banach space with $\dim X\ge2$, and let
$F:X\to X$ be a bijective step-isometry satisfying
\[
F(0)=0,
\qquad
\sup_{x\in X}\|F(x)-x\|<\infty.
\]
For every $v\ne0$, there exists a countable rational linear
subspace $A\subseteq X$ containing $v$, with $F(A)=A$, such
that $Z=\overline{A}^{\|\cdot\|}$ has dimension at least two and admits
a weak$^*$-compact set $K\subseteq B_{Z^*}$ satisfying
\begin{equation}\label{Knorming}
\|z\|=\sup_{\phi\in K}|\phi(z)|
\qquad(z\in Z)
\end{equation}
and
\begin{equation}\label{Korder}
\phi(Fx)\le\phi(Fy)
\end{equation}
for all $x,y\in A,\ \phi\in K$ with $\phi(x)<\phi(y)$.
\end{lemma}

\begin{proof}
Choose $w\in X$ linearly independent of $v$, and define
\[
A_0=\operatorname{span}_{\mathbb Q}\{v,w\},
\qquad
A_{n+1}
=\operatorname{span}_{\mathbb Q}
\bigl(A_n\cup F(A_n)\cup F^{-1}(A_n)\bigr).
\]
Let
\[
A=\bigcup_{n=0}^{\infty}A_n,
\qquad Z=\overline{A}^{\|\cdot\|}.
\]
Each $A_n$ is countable, and the sequence $(A_n)$ is
increasing. Hence $A$ is a countable rational linear
subspace containing $v$ and $w$. By construction,
\[
F(A)\subseteq A,
\qquad F^{-1}(A)\subseteq A,
\]
and therefore $F(A)=A$.

The closure \(Z=\overline{A}^{\|\cdot\|}\) is a real linear subspace of \(X\), since \(A\) is rational linear and \(\mathbb Q\) is dense in \(\mathbb R\). Thus \(Z\) is a separable Banach space  with countable dense subset \(A\). Since \(v,w\in Z\) are linearly independent, \(\dim Z\ge2\).

By Lemma~\ref{mazur-smoothness}, the G\^ateaux smooth
points of the norm of $Z$ are dense in $S_Z$. At each
such point $u$, let $j_Z(u)$ denote its unique unit
supporting functional, so that
\[
\|j_Z(u)\|=1,
\qquad j_Z(u)(u)=1.
\]
Define
\[
\Sigma
=\{j_Z(u):u\in S_Z\text{ is G\^ateaux smooth}\},
\qquad
K=\overline{\Sigma}^{\,w^*}.
\]
By the Banach-Alaoglu theorem,  \(K=\overline{\Sigma}^{\,w^*}\subseteq B_{Z^*}\) is weak$^*$-compact.

We now verify that $K$ determines the norm of $Z$. It is immediate that
\[
\sup_{\phi\in K}|\phi(z)|\le\|z\|
\qquad(z\in Z).
\]
because  $K\subseteq B_{Z^*}$ .
For the reverse inequality, fix $z\ne0$ and set
$e=\frac z{\|z\|}$. Choose smooth points $u_n\in S_Z$
with $\lim_{n\to\infty}u_n= e$, and write $\phi_n=j_Z(u_n)\in K$.
Then
\[
\begin{aligned}
\phi_n(e)=\phi_n(u_n)+\phi_n(e-u_n)
\ge1-\|\phi_n\|\,\|e-u_n\|=1-\|e-u_n\|.
\end{aligned}
\]
Since $\phi_n(e)\le1$, we obtain
\[
\lim_{n\to\infty}\phi_n(e)=1,
\qquad
\lim_{n\to\infty}\phi_n(z)=\|z\|.
\]
Consequently,
\[
\|z\|
=\lim_{n\to\infty}\phi_n(z)
\le\sup_{\phi\in K}|\phi(z)|
\le\|z\|,
\]
which establishes \eqref{Knorming}.

The restriction $F|_A:A\to A$ is a step-isometry
on a dense subset of $Z$, and
\[
\|F(x)-x\|_Z\le M
\qquad(x\in A).
\]
By Proposition~\ref{prop:order}, and the fact that $\Sigma$ is $w^*-$dense in $K$, we conclude
\[
\varphi(x)<\varphi(y)
\quad\Longrightarrow\quad
\varphi(Fx)\le\varphi(Fy)
\qquad(x,y\in A,\ \varphi\in K).
\]

\end{proof}

\begin{lemma}\label{fiber}
Let $A$ be a dense subset of a real Banach space $Z$  with $0\in A$, and let
$F:A\to A$ be a bijective step-isometry satisfying $F(0)=0$.
Suppose that $\phi\in Z^*\setminus\{0\}$ satisfies
\[
\phi(x)<\phi(y)\quad\Longrightarrow\quad
\phi(Fx)\le\phi(Fy)\qquad(x,y\in A).
\]
Then there exists a nondecreasing surjection $g_\phi:\phi(A)\to\phi(A)$ such that
\[
\phi(Fx)=g_\phi(\phi(x))\quad(x\in A),\qquad g_\phi(0)=0.
\]
Moreover, $g_\phi$ is continuous at 0 with respect to the relative topology of $\phi(A)$.
\end{lemma}

\begin{proof}
Let
\[
D:=\phi(A),\qquad
\rho:=\|\phi\|_{Z^*}>0,\qquad
H(x):=\phi(Fx)\quad(x\in A).
\]
The assumed order property gives
\begin{equation}\label{fiber-order}
\phi(x)<\phi(y)
\quad\Longrightarrow\quad
H(x)\le H(y)
\qquad(x,y\in A).
\end{equation}

We first prove that $H$ is constant on each level
set of $\phi$ in $A$.
For any $x\in A$, we have
\[
F\bigl(A\cap B_Z(x,1)\bigr)
=
A\cap B_Z(Fx,1).
\]
Consequently,
\[
\sup_{z\in A\cap B_Z(x,1)}H(z)
=\sup_{w\in A\cap B_Z(Fx,1)}\phi(w)=\phi(Fx)+\rho=H(x)+\rho.
\]

We claim that
\begin{equation}\label{fiberidentity}
H(x)+\rho
=
\sup\{H(z):z\in A,\ \phi(z)<\phi(x)+\rho\}.
\end{equation}
Indeed, if $z\in A\cap B_Z(x,1)$, then
\[
\phi(z)
\le \phi(x)+\rho\|z-x\|
<\phi(x)+\rho.
\]
Thus
\[
H(x)+\rho
\le
\sup\{H(z):z\in A,\ \phi(z)<\phi(x)+\rho\}.
\]
Conversely, suppose that $z\in A$ satisfies
$\phi(z)<\phi(x)+\rho$. There exists $w\in A\cap B_Z(x,1)$ such that
$\phi(z)<\phi(w)$ since
\(
\sup_{w\in A\cap B_Z(x,1)}\phi(w)
=
\phi(x)+\rho.
\)
By~\eqref{fiber-order},
\[
H(z)\le H(w)
\le
\sup_{u\in A\cap B_Z(x,1)}H(u)
=
H(x)+\rho.
\]
which establishes \eqref{fiberidentity}.

By \eqref{fiberidentity}, if \(\phi(x)=\phi(y)\), then \(H(x)=H(y).\)
Since \(F(A)=A\), we obtain a well-defined map \(g_\phi:D\to D\) given by
\[
g_\phi(\phi(x))=\phi(Fx)=H(x)\qquad(x\in A).
\]

If $a,b\in D$ satisfy $a<b$, choose $x,y\in A$
with $\phi(x)=a$ and $\phi(y)=b$.
By \eqref{fiber-order},  
\[
g_\phi(a)=H(x)\le H(y)=g_\phi(b).
\]
Thus $g_\phi$ is nondecreasing. Moreover,
\(
g_\phi(0)=\phi(F(0))=0.
\)

Let $b\in D$, and choose $y\in A$ such that $\phi(y)=b$.
Since $F(A)=A$, there exists $x\in A$ such that
$Fx=y$. Hence
\[
g_\phi(\phi(x))=\phi(Fx)=\phi(y)=b.
\]
Therefore $g_\phi(D)=D$.

Since $\phi\ne0$ and $A$ is dense in $Z$,
the set $D=\phi(A)$ is dense in $\mathbb R$.
Given $\varepsilon>0$, density and surjectivity
yields points $a_-,a_+\in D$ such that
\[
-\varepsilon<g_\phi(a_-)<0
<
g_\phi(a_+)<\varepsilon.
\]
Monotonicity together with $g_\phi(0)=0$ implies
\(
a_-<0<a_+.
\)
Let
\(
\delta:=\min\{-a_-,a_+\}>0.
\)
For any $a\in D$ with $|a|<\delta$, we have
\[
-\varepsilon
<
g_\phi(a_-)
\le g_\phi(a)
\le g_\phi(a_+)
<
\varepsilon.
\]
Thus \(g_\phi\) is continuous at 0 in \(D\).
\end{proof}

\begin{theorem}\label{continuity}
Every whole-space bijective step-isometry $T:X\to X$ is a homeomorphism.
\end{theorem}
\begin{proof}
By Lemma~\ref{lem:normalization}, there exists a surjective
linear isometry $U:X\to X$ such that
\[
F(x):=U^{-1}(Tx-T(0))
\]
is a whole-space bijective step-isometry satisfying
\[
F(0)=0,
\qquad
\|F(x)-x\|\le M
\qquad(x\in X)
\]
for some finite constant $M$.

 Fix $v\in X\setminus\{0\}$.
By Lemma~\ref{invariant-model}, there exist \(A\), \(Z=\overline{A}^{\|\cdot\|}\), and a weak$^*$-compact norming set \(K\subseteq B_{Z^*}\) such that
\[
F(A)=A,\qquad
\pm \frac vn\in A\quad(n\in\mathbb N),
\]
and, by~\eqref{Knorming},
\[
\|z\|=\sup_{\phi\in K}|\phi(z)|
\qquad(z\in Z).
\]

For each $n\in\mathbb N$, define
\[
h_n:K\to[0,\infty),
\qquad
h_n(\phi):=
\max\left\{
|\phi(F(\frac vn))|,
|\phi(F(-\frac vn))|
\right\}.
\]
Since $F(\pm \frac vn)\in A\subseteq Z$,   $h_n$ is weak$^*$-continuous on $K$.

For \(\phi\ne0\), by Lemma~\(\ref{fiber}\), with \(a=|\phi(v)|\), we have
\[
h_n(\phi)=
\max\left\{
g_\phi\!\left(\frac an\right),
-g_\phi\!\left(-\frac an\right)
\right\}.
\]
By monotonicity and \(g_\phi(0)=0\),
\[
0\le g_\phi\!\left(\frac{a}{n+1}\right)
\le g_\phi\!\left(\frac an\right),
\qquad
0\le-g_\phi\!\left(-\frac{a}{n+1}\right)
\le-g_\phi\!\left(-\frac an\right).
\]
 Therefore
\[
0\le h_{n+1}(\phi)\le h_n(\phi),
\qquad
\lim_{n\to\infty}h_n(\phi)=0.
\]

For $\varepsilon>0$, let
\[
V_n:=\{\phi\in K:h_n(\phi)<\varepsilon\}.
\]
This forms an increasing open cover of \(K\). By compactness, \(K=V_N\) for some \(N\). Hence,
\[
\max\left\{
\left\|F\!\left(\frac vn\right)\right\|,
\left\|F\!\left(-\frac vn\right)\right\|
\right\}
=\sup_{\phi\in K}h_n(\phi)
\le\varepsilon
\qquad(n\ge N).
\]
Thus
\[
\lim_{n\to\infty}\left\|F\!\left(\frac vn\right)\right\|=0
\qquad(v\in X),
\]
including \(v=0\). By lemmas~\(\ref{closedballs}\) and \(\ref{barrel}\), we know that \(F\) is continuous at 0.

For \(a\in X\), the map
\[
F_a(x)=F(a+x)-F(a)
\]
is a bijective step-isometry fixing 0, with
\[
\|F_a(x)-x\|
\le\|F(a+x)-(a+x)\|+\|F(a)-a\|
\le2M.
\]
The above argument shows that \(F_a\) is continuous at 0. Thus
\[
\lim_{x\to0}\|F(a+x)-F(a)\|
=\lim_{x\to0}\|F_a(x)\|=0.
\]
Since \(a\) was arbitrary, \(F\) is continuous on \(X\), and therefore so is \(T=UF+T(0)\). Applying the same continuity result to \(T^{-1}\), which is also a bijective step-isometry, proves that \(T\) is a homeomorphism.
\end{proof}

This answers Question~\ref{question3}. Completeness enters through the
Baire theorem on $X$. Separability is imposed only on the auxiliary space
$Z$, not on the original space $X$.

\section{Typical strongly non-Rado sets}\label{sec:typical}
Let $X$ be a real normed space, and let $S$ be countable and dense in $X$. The graph $G(S,p)$ retains each
unordered pair at distance less than one independently with probability
$p\in(0,1)$; all other pairs are nonedges. A graph is \emph{rigid} if its
only automorphism is the identity. We write $\Iso(S)$ for the group of
bijective isometries on $S$.

The graph distance is denoted by $d_G$.
The following distance-recovery observation is stated explicitly in
\cite[Section~1, discussion following Lemma~1]{Balogh}; see also
\cite[Lemma~3]{BBGLW}.
\begin{lemma}\label{graphdistance}
For $G=G(S,p)$, almost surely, simultaneously for all distinct $x,y\in S$
and integers $k\ge2$,
\(
d_G(x,y)\le k\) if and only if \(\norm{x-y}<k.
\)
\end{lemma}
The restriction $k\ge2$ is important. At distance less than one,
adjacency still has probability $p$, not probability one. We will use
rescaling to avoid assuming recovery of this first threshold.

\begin{lemma}[{\cite[Lemma~2]{Balogh}}]
\label{known-countable-step}
Let $S$ be a countable dense subset of a normed
space. If there are at most countably many bijective
step-isometries $S\to S$, then $S$ is strongly
non-Rado.
\end{lemma}

\begin{lemma}[{\cite[Theorem~12]{BroerePilsniak2015}}]
\label{known-random-edge}
Let $H$ be a countable graph with countable
automorphism group. Suppose that every nonidentity
automorphism of $H$ moves infinitely many edges.
For a fixed $p\in(0,1)$, colour each edge
independently black with probability $p$ and white
with probability $1-p$. Then, almost surely, the
identity is the only automorphism of $H$ preserving
this colouring.
\end{lemma}

The second quoted lemma is the fixed-probability special case of the
cited random edge-colouring theorem. Its colours will record retention
and deletion of potential edges.

\begin{proposition}\label{countablegroup}
Let $X$ be a real Banach space with $\dim X\ge2$
such that every bijective step-isometry of every
dense subset of $X$ is an isometry. Suppose that
$S\subseteq X$ is countable and dense, and that
$\operatorname{Iso}(S)$ is countable. Then:
\begin{enumerate}
\item
For any fixed $p,q\in(0,1)$, independent random
graphs $G(S,p)$ and $G(S,q)$ are almost surely
nonisomorphic.
\item
For each fixed $p\in(0,1)$, the graph $G(S,p)$
almost surely has no non-trivial automorphism.
\end{enumerate}
\end{proposition}
\begin{proof}
For each $r\in(0,1)$, Lemma~\ref{graphdistance} provides
a family $\Omega_r$ of probability one such that every
$G\in\Omega_r$ satisfies \(
d_G(x,y)\le k\) if and only if \(\norm{x-y}<k.
\)

Let $G_1\in\Omega_p$, $G_2\in\Omega_q$, and suppose that
$f:S\to S$ is an isomorphism from $G_1$ to $G_2$. Thus, for distinct $x,y\in S$ and every integer $k\ge2$,
\(\|x-y\|<k\) if and only if \(\|f(x)-f(y)\|<k\) since
\(
d_{G_1}(x,y)=d_{G_2}(f(x),f(y)).
\)

Put
\[
D:=\frac12 S
=\left\{\frac{x}{2}:x\in S\right\},
\qquad
\widetilde f\!\left(\frac{x}{2}\right)
:=\frac{f(x)}{2}.
\]
Then $\widetilde f:D\to D$ is bijective, and $D$ is dense in $X$.
For every positive integer $m$, we have $2m\ge2$, and hence
\[
\begin{aligned}
\left\|\frac{x}{2}-\frac{y}{2}\right\|<m
&\iff \|x-y\|<2m\\
&\iff \|f(x)-f(y)\|<2m\\
&\iff
\left\|
\widetilde f\!\left(\frac{x}{2}\right)
-\widetilde f\!\left(\frac{y}{2}\right)
\right\|<m.
\end{aligned}
\]
Therefore $\widetilde f$ is a bijective step-isometry.
Consequently,
\[
\begin{aligned}
\|f(x)-f(y)\|
=
2\left\|
\widetilde f\!\left(\frac{x}{2}\right)
-\widetilde f\!\left(\frac{y}{2}\right)
\right\|
=
2\left\|\frac{x}{2}-\frac{y}{2}\right\|
=\|x-y\|.
\end{aligned}
\]
Thus every such graph isomorphism belongs to
$\operatorname{Iso}(S)$. Taking $G_1=G_2$ gives the same
conclusion for automorphisms.

To prove \textup{(i)}, let $G_p=G(S,p)$ and
$G_q=G(S,q)$ be independent, and fix
$\sigma\in\operatorname{Iso}(S)$.
Choose infinitely many distinct allowed edges
\[
e_n=\{x_n,y_n\},
\qquad
0<\|x_n-y_n\|<1
\qquad(n\ge1).
\]
Such a sequence exists because $S$ is dense.
Since $\sigma$ is a bijective isometry, the pairs
\[
\sigma(e_n)=\{\sigma(x_n),\sigma(y_n)\}
\]
are also distinct allowed edges.

Define
\[
X_n:=\mathbf 1_{\{e_n\in E(G_p)\}},
\qquad
Y_n:=\mathbf 1_{\{\sigma(e_n)\in E(G_q)\}}.
\]
Independence of edge selection within each graph,
together with independence of the two graphs, shows
that the variables
\[
X_1,Y_1,X_2,Y_2,\ldots
\]
are independent, with
\[
\mathbb P(X_n=1)=p,
\qquad
\mathbb P(Y_n=1)=q.
\]
For $A_n:=\{X_n=Y_n\}$, we have
\[
\mathbb P(A_n)=
\mathbb P(X_n=1,Y_n=1)
+\mathbb P(X_n=0,Y_n=0)
=pq+(1-p)(1-q).
\]
Let \(\theta=pq+(1-p)(1-q)\), then
\(
1-\theta=p(1-q)+(1-p)q>0,
\)
we have $\theta<1$.

If $\sigma$ induces an isomorphism from $G_p$ to $G_q$,
then $A_n$ occurs for every $n$. Consequently, for every
$N\ge1$,
\[
\begin{aligned}
\mathbb P(\sigma:G_p\cong G_q)
\le
\mathbb P\left(\bigcap_{n=1}^{N}A_n\right)
=
\prod_{n=1}^{N}\mathbb P(A_n)
=\theta^N.
\end{aligned}
\]
Letting $N\to\infty$ gives
\[
\mathbb P(\sigma:G_p\cong G_q)=0.
\]

On the probability-one event
$\{G_p\in\Omega_p,\ G_q\in\Omega_q\}$, every isomorphism
belongs to $\operatorname{Iso}(S)$. Hence
\[
\begin{aligned}
\mathbb P(G_p\cong G_q)
\le
\sum_{\sigma\in\operatorname{Iso}(S)}
\mathbb P(\sigma:G_p\cong G_q)
=0,
\end{aligned}
\]
where countability of $\operatorname{Iso}(S)$ is used.
This proves \textup{(i)}.

To prove \textup{(ii)}, fix
\[
\sigma\in
\operatorname{Iso}(S)\setminus\{\mathrm{id}_S\}.
\]
Choose $x\in S$ such that $\sigma(x)\ne x$.
Since $S\cap B(x,1)$ is infinite, we can choose distinct
points
\[
y_n\in
\bigl(S\cap B(x,1)\bigr)
\setminus\{x,\sigma(x),\sigma^{-1}(x)\}.
\]
Put
\[
e_n:=\{x,y_n\},
\qquad
e_n':=\{\sigma(x),\sigma(y_n)\}.
\]
Both are allowed edges, since
\[
0<\|x-y_n\|
=\|\sigma(x)-\sigma(y_n)\|<1.
\]

The edges $e_n$ are pairwise distinct, and injectivity
of $\sigma$ implies that the edges $e_n'$ are pairwise
distinct. Furthermore,
\[
x\in e_n,
\qquad
x\notin e_m'
\qquad(n,m\ge1),
\]
because
\[
\sigma(x)\ne x,
\qquad
y_m\ne\sigma^{-1}(x).
\]
It follows that
\[
e_n\ne e_m'
\qquad(n,m\ge1).
\]
Thus all edges in the collection
\[
e_1,e_1',e_2,e_2',\ldots
\]
are distinct.

For $G=G(S,p)$, define
\[
U_n:=\mathbf 1_{\{e_n\in E(G)\}},
\qquad
V_n:=\mathbf 1_{\{e_n'\in E(G)\}}.
\]
These are independent Bernoulli variables with parameter
$p$. Therefore the events
\[
B_n:=\{U_n=V_n\}
\]
are independent and satisfy
\[
\mathbb P(B_n)
=p^2+(1-p)^2
=:\eta<1.
\]
Here independence follows from distinctness of the
edges; disjointness of their vertex sets is not required.

If $\sigma\in\operatorname{Aut}(G)$, then $B_n$ occurs
for every $n$. Hence, for every $N\ge1$,
\[
\begin{aligned}
\mathbb P(\sigma\in\operatorname{Aut}(G))
\le
\mathbb P\left(\bigcap_{n=1}^{N}B_n\right)
=
\prod_{n=1}^{N}\mathbb P(B_n)
=\eta^N.
\end{aligned}
\]
Letting $N\to\infty$, we obtain
\[
\mathbb P(\sigma\in\operatorname{Aut}(G))=0.
\]

On $\Omega_p$, every automorphism of $G$ belongs to
$\operatorname{Iso}(S)$. Countability now gives
\[
\begin{aligned}
\mathbb P\bigl(
\operatorname{Aut}(G)\ne\{\mathrm{id}_S\}
\bigr)
\le
\sum_{\substack{
\sigma\in\operatorname{Iso}(S)\\
\sigma\ne\mathrm{id}_S}}
\mathbb P(\sigma\in\operatorname{Aut}(G))
=0.
\end{aligned}
\]
This proves \textup{(ii)}.
\end{proof}

The following is the standard finite-base argument; compare the discussion
following Theorem~5.3 in~\cite{Lehner2013}. The short proof also records its
metric consequence.
\begin{lemma}\label{finitebase}
Let $S\subseteq X$ be countable. If a finite set $B=\{b_1,\ldots,b_m\}\subset S$ is fixed pointwise only by the identity in $\Iso(S)$, then $\Iso(S)$ is countable. In particular, this holds if
\[
x\mapsto(\norm{x-b_1},\ldots,\norm{x-b_m})
\]
is injective on $S$. One point $b_1$ with distinct distances to all other points suffices.
\end{lemma}
\begin{proof}
Define
\[
\Phi:\operatorname{Iso}(S)\longrightarrow S^m,
\qquad
\Phi(\sigma)
=\bigl(\sigma(b_1),\ldots,\sigma(b_m)\bigr).
\]
If $\Phi(\sigma)=\Phi(\tau)$, then
\[
(\sigma^{-1}\circ\tau)(b_i)
=\sigma^{-1}(\sigma(b_i))
=b_i
\qquad(1\le i\le m).
\]
By hypothesis, $\sigma^{-1}\circ\tau=\mathrm{id}_S$,
so $\sigma=\tau$. Thus $\Phi$ is injective.
Since $S$ is countable and $m$ is finite, $S^m$ is
countable, and hence so is $\operatorname{Iso}(S)$.

Suppose now that
\[
D_B:S\longrightarrow\mathbb R^m,
\qquad
D_B(x)=\bigl(\|x-b_1\|,\ldots,\|x-b_m\|\bigr)
\]
is injective. Let $\sigma\in\operatorname{Iso}(S)$
fix every point of $B$. For every $x\in S$ and
$1\le i\le m$,
\[
\|\sigma(x)-b_i\|
=\|\sigma(x)-\sigma(b_i)\|
=\|x-b_i\|.
\]
Therefore
\[
D_B(\sigma(x))=D_B(x).
\]
Injectivity of $D_B$ gives $\sigma(x)=x$ for every
$x\in S$. Thus only the identity fixes $B$ pointwise,
and the preceding argument applies.

Finally, suppose that the distances
\[
\{\|x-b_1\|:x\in S\setminus\{b_1\}\}
\]
are pairwise distinct. They are positive, whereas
$\|b_1-b_1\|=0$. Hence
\[
x\longmapsto\|x-b_1\|
\]
is injective on $S$, and the conclusion follows
by taking $B=\{b_1\}$.
\end{proof}

We use two precise formulations of typicality. Category statements concern
sequences in the product space $X^{\N}$. Sampling statements concern the
product law $\mu^{\N}$ of an explicitly specified Borel probability.
A probability has \emph{full support} if every nonempty open set has positive
measure. In the null-sphere hypothesis, a sphere means
$\{x:\norm{x-c}=r\}$ with $c\in X$ and $r\ge0$.

\begin{theorem}\label{category}
For every separable Banach space $X$, with $X^{\N}$ given the product topology, there is a dense $G_\delta$ subset of $X^{\N}$ consisting of sequences with distinct terms, dense range, and
\begin{equation}\label{noisosceles}
\norm{x_i-x_k}\ne\norm{x_j-x_k}\qquad(i,j,k\text{ pairwise distinct}).
\end{equation}
If every step-isometry on every dense self-domain of $X$ is an isometry, the resulting sets are strongly non-Rado and for each fixed $p\in(0,1)$ their graphs $G(S,p)$ are almost surely rigid.
\end{theorem}
\begin{proof}
Equip $X^{\mathbb N}$ with the product topology.
It is completely metrizable, for example by
\[
d(\mathbf x,\mathbf y)
=
\sum_{n=1}^{\infty}
2^{-n}\min\{1,\|x_n-y_n\|\},
\]
and is therefore a Baire space.

Since $X$ is separable, choose a countable basis
$(U_\ell)_{\ell\ge1}$ of nonempty open subsets of $X$.
For each $\ell$, put
\[
\mathcal V_\ell
:=
\{\mathbf x\in X^{\mathbb N}:
x_n\in U_\ell\text{ for some }n\}.
\]
This set is open. It is also dense, because a basic
open cylinder restricts only finitely many
coordinates, leaving another coordinate free
to take a value in $U_\ell$.
Consequently, every sequence in
$\bigcap_{\ell\ge1}\mathcal V_\ell$ has dense range.

For $i<j$, define
\[
\mathcal D_{ij}
:=
\{\mathbf x\in X^{\mathbb N}:x_i=x_j\}.
\]
This set is closed and has empty interior, since
a coordinate can be perturbed within any prescribed
open neighbourhood to avoid a fixed point.

For pairwise distinct $i,j,k$, define
\[
\mathcal E_{ijk}
:=
\{\mathbf x\in X^{\mathbb N}:
\|x_i-x_k\|=\|x_j-x_k\|\}.
\]
Continuity of the norm shows that
$\mathcal E_{ijk}$ is closed.
We claim that it has empty interior.
Indeed, in any nonempty basic open cylinder,
fix all coordinates except $x_i$.
The admissible values of $x_i$ form a nonempty
open subset of $X$, whereas the forbidden values
lie on the sphere
\[
\{u\in X:\|u-x_k\|=\|x_j-x_k\|\}.
\]
Every norm sphere, including a radius-zero sphere,
has empty interior. Hence an admissible value
outside this sphere can be chosen.
Thus $\mathcal E_{ijk}$ is nowhere dense.

It follows from the Baire category theorem that
\[
\mathcal R
:=
\left(\bigcap_{\ell\ge1}\mathcal V_\ell\right)
\cap
\left(\bigcap_{i<j}
(X^{\mathbb N}\setminus\mathcal D_{ij})\right)
\cap
\left(
\bigcap_{\substack{i,j,k\\
                   \text{pairwise distinct}}}
(X^{\mathbb N}\setminus\mathcal E_{ijk})
\right)
\]
is a dense $G_\delta$ subset of $X^{\mathbb N}$.
Every sequence in $\mathcal R$ has distinct terms,
dense range, and satisfies~\eqref{noisosceles}.

Now suppose that every bijective step-isometry
on every dense self-domain of $X$ is an isometry.
Fix $\mathbf x=(x_n)\in\mathcal R$, and put
\[
S:=\{x_n:n\in\mathbb N\}.
\]
We claim that
\[
S\longrightarrow[0,\infty),
\qquad
s\longmapsto\|s-x_1\|
\]
is injective.
The point $x_1$ is the unique point at distance
zero from itself. For distinct $x_i,x_j\ne x_1$,
condition~\eqref{noisosceles} gives
\[
\|x_i-x_1\|\ne\|x_j-x_1\|.
\]
The claim follows.

Lemma~\ref{finitebase}, applied to the single
anchor $x_1$, therefore shows that
$\operatorname{Iso}(S)$ is countable.
Proposition~\ref{countablegroup} now implies that
$S$ is strongly non-Rado and that, for each fixed
$p\in(0,1)$, $G(S,p)$ almost surely has no
nonidentity automorphism.
\end{proof}

The no-isosceles condition compares two distances with a common endpoint.
It gives a trivial single-point stabilizer and hence a countable isometry
group; no assertion that the whole isometry group is trivial is needed.

\begin{theorem}\label{sampling}
Let $\mu$ be a full-support Borel probability on a separable Banach space $X$, and assume every sphere, including radius-zero spheres, has $\mu$-measure zero. An independent sequence with law $\mu$ almost surely has distinct terms, dense range, and no repeated distances between distinct unordered pairs. If every step-isometry on every dense self-domain of $X$ is an isometry, the sampled range is almost surely strongly non-Rado and for each fixed $p\in(0,1)$ its graphs $G(S,p)$ are almost surely rigid.
\end{theorem}
\begin{proof}
Let $(Y_n)_{n\ge1}$ be independent random variables
with common law $\mu$.

Since radius-zero spheres have measure zero,
$\mu$ has no atoms. Thus, for $i\ne j$,
\[
\mathbb P(Y_i=Y_j)
=
\int_X\mu(\{y\})\,d\mu(y)
=
0.
\]
Taking the union over the countably many pairs
$i<j$ shows that the terms are almost surely
distinct.

Next, choose a countable basis $(U_\ell)_{\ell\ge1}$
of nonempty open subsets of $X$.
Full support gives $\mu(U_\ell)>0$.
By independence,
\[
\begin{aligned}
\mathbb P(Y_n\notin U_\ell\text{ for every }n)
&=
\lim_{N\to\infty}
\mathbb P(Y_1,\ldots,Y_N\notin U_\ell)\\
&=
\lim_{N\to\infty}(1-\mu(U_\ell))^N\\
&=0.
\end{aligned}
\]
Since the basis is countable, almost surely the
sampled range meets every $U_\ell$ and is therefore
dense in $X$.

We now exclude repeated distances.
Fix two distinct unordered index pairs
\[
\{i,j\}\ne\{k,\ell\},
\qquad i\ne j,\quad k\ne\ell.
\]
After relabelling their endpoints, we may assume
that $i\notin\{k,\ell\}$. Consider the event
\[
E:=\{\|Y_i-Y_j\|=\|Y_k-Y_\ell\|\}.
\]
Condition on all the other random variables
appearing in this equality. The remaining variable
$Y_i$ is independent of them and still has law
$\mu$. Hence
\[
\begin{aligned}
\mathbb P(E\mid Y_j,Y_k,Y_\ell)
&=
\mu\bigl(
\{z\in X:\|z-Y_j\|=\|Y_k-Y_\ell\|\}
\bigr)\\
&=0
\qquad\text{almost surely}.
\end{aligned}
\]
The last equality follows from the null-sphere
assumption, including the possible case of
radius zero. Therefore $\mathbb P(E)=0$.

There are only countably many choices of two
distinct unordered index pairs. Consequently,
almost surely,
\[
\|Y_i-Y_j\|\ne\|Y_k-Y_\ell\|
\quad\text{whenever}\quad
\{i,j\}\ne\{k,\ell\},
\quad i\ne j,\quad k\ne\ell.
\]
Together with the preceding conclusions, this
gives a probability-one event on which the terms
are distinct, their range is dense, and all
distances between distinct unordered pairs
are different.

Now assume the stated rigidity property of $X$,
and fix a sample in this probability-one event.
Put
\[
S:=\{Y_n:n\in\mathbb N\}.
\]
The map
\[
S\to[0,\infty),
\qquad s\mapsto\|s-Y_1\|
\]
is injective. Indeed, $Y_1$ is the unique point
at distance zero from itself, and all other
distances from $Y_1$ are distinct.
Lemma~\ref{finitebase} therefore implies that
$\operatorname{Iso}(S)$ is countable.

Proposition~\ref{countablegroup} now shows that
$S$ is strongly non-Rado and that, for each fixed
$p\in(0,1)$, $G(S,p)$ almost surely has no
nonidentity automorphism.
\end{proof}

The probability-one event in this theorem concerns the sampled sequence.
For each range $S$ in that event and each fixed edge parameter, the graph
conclusion concerns subsequent independent edge sampling on $S$.

The following classical fact is used without proof.
\begin{lemma}[{\cite[Corollary~2.1]{Borell1974}}]
\label{gaussian-log-concavity}
Every Radon Gaussian probability measure $\mu$
on a real Banach space is log-concave. In particular,
for open convex sets $A,B$ and $0<t<1$,
\[
\mu((1-t)A+tB)
\ge
\mu(A)^{1-t}\mu(B)^t.
\]
\end{lemma}
\begin{proposition}\label{Gaussian}
Let $X$ be a separable real Banach space. Choose unit vectors $(u_n)$ with dense linear span and independent standard normal real variables $(Z_n)$. The series
\[
Y=\sum_{n=1}^\infty2^{-n}Z_nu_n
\]
converges almost surely. Its distribution is a full-support Gaussian Borel probability assigning zero mass to every norm sphere. Strict convexity is not required.
\end{proposition}
\begin{proof}
Since $\|u_n\|=1$, Tonelli's theorem gives
\[
\begin{aligned}
\mathbb E\left[
\sum_{n=1}^{\infty}\|2^{-n}Z_nu_n\|
\right]
&=
\sum_{n=1}^{\infty}2^{-n}\mathbb E|Z_n|\\
&=
\sqrt{\frac2\pi}
\sum_{n=1}^{\infty}2^{-n}
<\infty.
\end{aligned}
\]
Hence the series converges absolutely almost surely.
Completeness of $X$ gives an $X$-valued limit $Y$.
The event of convergence is measurable. Define $Y=0$
on its null complement. The resulting $Y$ is a Borel
random variable, being the limit of the measurable
partial sums on the convergence event.

For every $\phi\in X^*$,
\[
\phi(Y)
=
\sum_{n=1}^{\infty}2^{-n}Z_n\phi(u_n).
\]
The coefficients are square summable, so
\[
\phi(Y)\sim
N\left(0,\sum_{n=1}^{\infty}
4^{-n}\phi(u_n)^2\right).
\]
Thus the distribution $\mu$ of $Y$ is Gaussian.
If $\phi\ne0$, density of the linear span of
$(u_n)$ ensures that $\phi(u_n)\ne0$ for some $n$.
Therefore the displayed variance is positive.

We next prove full support.
Fix $x\in X$ and $\varepsilon>0$.
Choose a finite linear combination
\[
v=\sum_{n=1}^{K}\alpha_nu_n,
\qquad
\|x-v\|<\varepsilon/3.
\]
Write
\[
H_N:=\sum_{n=1}^{N}2^{-n}Z_nu_n,
\qquad
R_N:=Y-H_N.
\]
Choose $N\ge K$ sufficiently large that
\[
\mathbb E\|R_N\|
\le
\sqrt{\frac2\pi}\sum_{n>N}2^{-n}
<
\varepsilon/6.
\]
Markov's inequality yields
\[
\mathbb P(\|R_N\|<\varepsilon/3)>\frac12.
\]

Set $\alpha_n=0$ for $K<n\le N$.
The event
\[
A_N:=
\bigcap_{n=1}^{N}
\left\{
|2^{-n}Z_n-\alpha_n|
<
\frac{\varepsilon}{3N}
\right\}
\]
has positive probability and implies
$\|H_N-v\|<\varepsilon/3$.
Since $A_N$ and $R_N$ are independent,
\[
\begin{aligned}
\mu(B(x,\varepsilon))
&\ge
\mathbb P\bigl(
A_N\cap\{\|R_N\|<\varepsilon/3\}
\bigr)\\
&=
\mathbb P(A_N)
\mathbb P(\|R_N\|<\varepsilon/3)
>0.
\end{aligned}
\]
Hence $\mu$ has full support.

It remains to prove that spheres are $\mu$-null.
As a Borel probability on the Polish space $X$,
$\mu$ is Radon. Fix $c\in X$ and put
\[
f_c(r):=\mu(B(c,r))
\qquad(r>0).
\]
Full support gives $f_c(r)>0$.
By Lemma~\ref{gaussian-log-concavity} and
\[
(1-t)B(c,r)+tB(c,s)
=
B(c,(1-t)r+ts),
\]
we obtain
\[
f_c((1-t)r+ts)
\ge
f_c(r)^{1-t}f_c(s)^t
\qquad(r,s>0,\ 0<t<1).
\]
Thus $\log f_c$ is a finite concave function on
$(0,\infty)$ and is therefore continuous.
Continuity from above now gives, for every $r>0$,
\[
\begin{aligned}
\mu(\overline B(c,r))
&=
\lim_{m\to\infty}\mu(B(c,r+1/m))\\
&=
\lim_{m\to\infty}f_c(r+1/m)\\
&=
f_c(r)
=
\mu(B(c,r)).
\end{aligned}
\]
Consequently,
\[
\mu\{x\in X:\|x-c\|=r\}=0.
\]

For radius zero, choose a nonzero $\phi\in X^*$.
The nondegenerate Gaussian variable $\phi(Y)$
has no atoms, so
\[
\mu(\{c\})
\le
\mathbb P(\phi(Y)=\phi(c))
=
0.
\]
This completes the proof.
\end{proof}

Only convexity of norm balls is used in the last step; strict convexity
is not required for this sampling model.

We can now apply the extension theorem to the two general constructions.
\begin{theorem}\label{typicalmain}
Let $X$ be strictly convex with norm-separable dual.
\begin{enumerate}
\item There is a dense $G_\delta$ subset of $X^{\N}$ consisting of sequences with distinct terms and dense range $S$ such that $S$ is strongly non-Rado.
\item If $\mu$ is a full-support Borel probability assigning zero mass to every sphere, independent sampling with law $\mu$ almost surely produces a countable dense strongly non-Rado set $S$.
\end{enumerate}
For every such $S$ and each fixed $p\in(0,1)$, $G(S,p)$ is almost surely rigid. For any fixed $p,q\in(0,1)$, independent $G(S,p)$ and $G(S,q)$ are almost surely nonisomorphic.
\end{theorem}

\begin{proof}
Since $X^*$ is norm-separable, it follows that  $X$ is separable and Asplund, moreover, the Fr\'echet smooth points of the
given norm are norm dense in $S_X$.
By Theorem~\ref{yangetuyantuo}, every bijective step-isometry on every dense
self-domain of $X$ is the restriction of a surjective affine isometry,
and hence is an isometry.

Theorem~\ref{category} therefore proves \textup{(i)}.
For \textup{(ii)}, apply Theorem~\ref{sampling} to the stated probability
$\mu$. In both constructions, distances from an anchor distinguish the
points of $S$, so Lemma~\ref{finitebase} makes $\Iso(S)$ countable.
Proposition~\ref{countablegroup} gives almost sure graph rigidity and
nonisomorphism for any two fixed edge parameters.
Proposition~\ref{Gaussian} supplies probabilities satisfying the sampling
assumptions.
\end{proof}
Theorem~\ref{typicalmain} proves Conjecture~\ref{conj:typical} in the two
specified senses of typicality.

\begin{corollary}\label{ell1typical}
The conclusions of Theorem~\ref{typicalmain} hold for standard real $\ell_1$. Thus Question~\ref{q:l1} has the strongly non-Rado answer in both the category and the stated sampling formulations.
\end{corollary}
\begin{proof}
Standard real $\ell_1$ is separable.
Let $S\subseteq\ell_1$ be dense and let
$T:S\to S$ be a bijective step-isometry.
Viewed as a map into $\ell_1$, $T$ has dense image
and is therefore almost surjective.
By Theorem \ref{main-lone-finite}, $T$ extends
to a surjective affine isometry of $\ell_1$.
Thus every bijective step-isometry on every dense
self-domain of $\ell_1$ is an isometry.

Theorems~\ref{category} and~\ref{sampling}, together
with Proposition~\ref{countablegroup}, now give
the asserted category and sampling conclusions,
including almost sure rigidity and nonisomorphism
for different fixed edge parameters.
Proposition~\ref{Gaussian} supplies full-support
sampling measures assigning zero mass to every sphere.
\end{proof}

For a concrete model in $\ell_1$, take $u_n=e_n$ in
Proposition~\ref{Gaussian}; then
\[
Y=(2^{-n}Z_n)_{n\ge1}\in\ell_1\quad\text{almost surely}.
\]
The argument for $\ell_1$ uses its representation theorem, rather than
the separable-dual hypothesis of Theorem~\ref{extension}.

\begin{remark}
These typical-set conclusions allow the exceptional Rado sets in $\ell_1$
and $\ell_2$ constructed in~\cite[Theorem~5]{Balogh}.
Full support alone does not suffice for the sampling conclusion:
a full-support atomic probability carried by a countable dense Rado set
almost surely samples every point of that set, so its sampled range is
Rado. The null-sphere assumption excludes this example.

For each $S$ obtained above and any prescribed countable sequence of fixed edge parameters, the
corresponding independent graphs are simultaneously rigid and pairwise
nonisomorphic almost surely, by taking a countable intersection of the
probability-one events. No common event for an uncountable family of
parameters is asserted.
\end{remark}

\end{document}